\documentclass[]{imsart}
\RequirePackage{amsthm,amsmath,amsfonts,amssymb}
\RequirePackage[numbers,sort&compress]{natbib}
\RequirePackage[colorlinks,citecolor=blue,urlcolor=blue]{hyperref}
\RequirePackage{graphicx}

\usepackage[margin=1.2in]{geometry}  

\startlocaldefs

\newtheorem{thm}{Theorem}[]
\newtheorem{lemma}[thm]{Lemma}
\newtheorem{remark}[thm]{Remark}

\renewcommand{\cite}{\citet*}
\newcommand{\bI}{{\bf I}}
\newcommand{\bX}{{\bf X}}
\newcommand{\bY}{{\bf Y}}
\newcommand{\bA}{{\bf A}}
\newcommand{\bM}{{\bf M}}
\newcommand{\bB}{{\bf B}}
\newcommand{\bC}{{\bf C}}
\newcommand{\bS}{{\bf S}}
\newcommand{\bH}{{\bf H}}
\newcommand{\bbH}{{\mathbb H}}
\newcommand{\bbM}{{\mathbb M}}
\newcommand{\bx}{{\bf x}}
\newcommand{\bu}{{\bf u}}

\newcommand{\bSig}{{\boldsymbol\Sigma}}
\newcommand{\bGa}{{\boldsymbol\Gamma}}

\renewcommand{\(}{\left(}
\renewcommand{\)}{\right)}
\newcommand{\re}{{\rm E}}
\newcommand{\rtr}{{\rm tr}}

\allowdisplaybreaks

\endlocaldefs

\begin{document}

\begin{frontmatter}
\title{Liberating dimension and  spectral norm: A universal approach to spectral properties of sample covariance matrices}
\runtitle{Harmonic approach to spectral properties of sample covariance matrices}

\begin{aug}

\author[A]{\fnms{Yanqing} \snm{ Yin}
			\ead[label=e3]{yinyq@cqu.edu.cn}}

\runauthor{Yanqing Yin}

\address[A]{School of Mathematics and Statistics,
			Chongqing University\printead[presep={,\ }]{e3}}

\end{aug}

\begin{abstract}
In this paper, our objective is to present a constraining principle governing the spectral properties of the sample covariance matrix. This principle exhibits harmonious behavior across diverse limiting frameworks, eliminating the need for constraints on the rates of dimension $p$ and sample size $n$, as long as they both tend to infinity. We accomplish this by employing a suitable normalization technique on the original sample covariance matrix.  Following this, we establish a harmonic central limit theorem for linear spectral statistics within this expansive framework. This achievement effectively eliminates the necessity for a bounded spectral norm on the population covariance matrix and relaxes constraints on the rates of dimension $p$ and sample size $n$, thereby significantly broadening the applicability of these results in the field of high-dimensional statistics. We illustrate the power of the established results by considering the test for covariance structure under high dimensionality, freeing both $p$ and $n$.
  
\end{abstract}

\begin{keyword}[class=MSC]
\kwd[Primary ]{62H15}
\kwd{62B20}
\kwd[; secondary ]{62D10}
\end{keyword}

\begin{keyword}
\kwd{ultra-high dimension
  covariance matrix; central limit theorem; limiting spectral distribution; linear spectral statistics; M-P law}
\end{keyword}

\end{frontmatter}

\section{Introduction and motivation}
The surge in high-dimensional statistics has become a focal point in recent years, spurred by the rapid advancements in data science. As datasets continue to evolve with higher dimensions, classical multivariate statistical tools, initially formulated for fixed dimensions, grapple with substantial challenges. The inadequacy of these traditional frameworks in handling the complexities of modern datasets has prompted a shift in focus towards the development of statistical methodologies capable of withstanding the expanding dimensions of variables. Random matrix theory, with its emphasis on exploring the spectral properties of matrices with random entries, has emerged as a crucial contributor to addressing these challenges.

Commencing with the analysis of the spectral properties of the Wishart matrix by Wishart in 1928 \cite{Wishart28G}, researchers have undertaken considerable efforts to extend the applicability of random matrix theory to divergent datasets in practical scenarios. As a crucial step in this endeavor, attention is often directed towards understanding the first-order limiting behaviors, such as the spreading of spectra. The Stieltjes transform method proves to be an invaluable tool in this context.
 Consider a Hermitian matrix $\mathbf{A}$ with dimensions $p \times p$ and eigenvalues denoted by $\lambda_j$, where $j = 1, 2, \ldots, p$. The empirical spectral distribution (ESD) of $\mathbf{A}$ is defined as
$F^{\mathbf{A}}(x) = \frac{1}{p}\sum_{j=1}^p I(\lambda_j \leq x),$
where $I(D)$ represents the indicator function of an event $D$. The Stieltjes transform of $F^{\mathbf{A}}(x)$ is then given by
$ m(z) = \int_{-\infty}^{+\infty} \frac{1}{x - z} dF^{\mathbf{A}}(x), $
with $z = u + \upsilon i \in \mathbb{C}^+$.
Now, let's consider the sample covariance matrix given by $$\mathbf{S}_n^0 = \frac{1}{n}\sum_{i=1}^n \mathbf{x}_i\mathbf{x}_i^* = \frac{1}{n}\mathbf{X}_n\mathbf{X}_n^*,$$ where `$*$' denotes the conjugate transpose and $\mathbf{X}_n=(\mathbf{x}_1,\mathbf{x}_2,\ldots,\mathbf{x}_n)$ is the observation data matrix.
The distinguished Marcenko-Pastur (M-P) law initially established in \cite{marchenko1967distribution} and later extended in \cite{Wachter78S, Silverstein95S}, asserts convergence under the following key conditions:
\begin{itemize}
	\item The independent component structure (ICS) condition: $\mathbf{X}_n=\boldsymbol{\Sigma}_n^{1/2}\mathbf{Y}_n$, where $\boldsymbol{\Sigma}_n^{1/2}$ represents the $p\times p$ Hermitian square root of the population covariance matrix $\boldsymbol{\Sigma}_n$ and $\mathbf{Y}_n=(\mathbf{y}_1,\mathbf{y}_2,\ldots,\mathbf{y}_n)$ is a $p\times n$ matrix with i.i.d standard complex random variables. Also, $F^{\boldsymbol{\Sigma}_n}\xrightarrow{d}H$, and the sequence $(\boldsymbol{\Sigma}_n)_n$ is bounded in spectral norm.
	\item The high dimensional framework: $p/n\to c\in(0,\infty)$.
\end{itemize}
Under these conditions, it can be confidently established that the ESD $F^{\mathbf{S}_n}$ of the sample covariance matrix $\mathbf{S}_n^0=\frac{1}{n}\boldsymbol{\Sigma}_n^{1/2}\mathbf{Y}_n\mathbf{Y}_n^*\boldsymbol{\Sigma}_n^{1/2}$ almost surely converges weakly to a nonrandom probability density function $F^0$ as $n\to \infty$. For each $z\in\mathbb{C}^+$, $m^0(z)=m_{F^0}^0(z)$ is a solution to the equation
\begin{align}\label{al3}
m^0(z)=\int\frac1{t(1-c-czm^0(z))-z}dH(t),
\end{align}
and this solution is unique in the set $\left\{m(z)\in\mathbb{C}^+: -(1-c)/z+cm^0(z)\in\mathbb{C}^+\right\}$.
This law reveals that, under the limiting framework where dimension divergent at the same rate with the sample size, although the sample eigenvalue is no longer consistent estimate the  corresponding population one, the distance between them is truly constant order. And one shall thus make statistical inference for population covariance matrix using the sample covariance matrix by noticing their unique relationship. This high dimensional framework thus gain more and more attention and there have been huge fruitful work under this framework, to name but a few \cite{BaiS10S,Bai99M,TracyW09D,YaoB15L,BloemendalK16P,Bao2022,Johnstone01D}.
 
However, in many practical statistical situations, such as in biostatistics, researchers encounter datasets with ultrahigh dimensions, where $p/n\to \infty$. An immediate concern arises: how do the theoretical results established under the classical high-dimensional framework change in the context of ultrahigh dimensions? Insight from the investigation in \cite{BaiY88C} reveals that under such conditions, the limiting spectral distribution (LSD) of the normalized sample covariance matrix $\frac{1}{\sqrt{n p}}\left(\mathbf{Y}_n^{T} \mathbf{Y}_n-p \mathbf{I}_n\right)$ converges to the semicircle law, which is the limiting law of the standard Wigner matrix. Notably, both $\mathbf{Y}_n^{T} \mathbf{Y}_n$ and $\mathbf{Y}_n\mathbf{Y}_n^{T}$ share the same non-zero eigenvalues.
This observation implies that when $p/n$ is large, apart from the $p-n$ zero eigenvalues, the distance between the spectra of the classical sample covariance matrix $\frac1n\mathbf{Y}_n\mathbf{Y}_n^{T}$ and the population covariance matrix $\mathbf{I}_p$ is on the order of $p/n$. This consistency aligns with the high-dimensional framework, where asymptotically, $p/n$ remains of a constant order. Furthermore, this result is extended in \cite{WangP14L,ChenP15C} to more general population covariance matrix cases. Additionally, central limit theorems (CLT) for linear spectral statistics (LSS) for these models are established in \cite{ChenP15C, Qiu2023}. In a recent study, \cite{Ding2023} delves into the model $\frac{1}{\sqrt{np}}_n\mathbf{Y}_n^*\boldsymbol{\Sigma}_n\mathbf{Y}$ directly, omitting the centralized matrix term. It is noteworthy that the LSD of such a model does not exist when $p/n$ diverges. Additionally, they make the assumption of an exact asymptotic framework where $p$ is of the order of $n^{\alpha}$ for some $\alpha\in(0,\infty)$ and the assumption that $\boldsymbol{\Sigma}_n$ is diagonal with entries bounded from both below and above.

As a pivotal point, in practical scenarios, both the sample size $n$ and the dimension $p$ are treated as known constants. The focus of asymptotic frameworks is geared towards achieving valuable information of the spectral properties of the sample covariance matrix under consideration in finite sample scenarios. In light of this, we aim to derive harmonic results, removing restrictions on the rates of divergence for both $n$ and $p$. To accomplish these, drawing inspiration from the Bai-Yin law, we propose a renormalization of the original sample covariance matrix $\mathbf{X}_n\mathbf{X}_n^{*}$ using the quantity $\left(\sqrt{p}+\sqrt{n}\right)^2\|\boldsymbol{\Sigma}_n\|$. It is noteworthy that under such normalization, the LSD of the sample covariance matrix always exists as long as both $p$ and $n$ tend to infinity. The specific LSD is uniquely determined by the rates of divergence for $p$ and $n.$ When they are of significantly different orders, the LSD degenerates into trivial cases. We categorize the scenarios as follows: as $p,n\to \infty$,
\begin{itemize}
	\item $p/n\to c\in (0,+\infty)$ is termed the comparable scenario.
	\item $p/n\to 0$ represents the large $n$ scenario.
	\item $p/n\to +\infty$ characterizes the large $p$ scenario.
\end{itemize}

A careful examination of the normalization parameter  reveals an additional noteworthy advantage. This normalization method offers a potential solution to the persistent challenge of an unbounded spectral norm in the population covariance matrix-a common occurrence in practical applications. In general, the application of random matrix theory in high-dimensional statistics necessitated the assurance that the population spectrum remains non-divergent. This assumption underlies many theoretical results, particularly those concerning second-order spectral properties such as the CLT for LSSs.

In situations where a spiked model is assumed, featuring finite but divergent spiked population eigenvalues, \cite{liu2023clt} has provided a partial solution. This is achieved by considering the joint distribution of the spiked eigenvalues and the CLT for LSS of the bounded, non-spiked population. For the general case, the findings in \cite{Yin2022} offer valuable insights. Specifically, \cite{Yin2022} identifies certain spectral statistics, showing that the convergence rate of CLT for LSS of the sample covariance matrix can be significantly influenced by the rate of the spectral norm under scenarios with a divergent spectrum. 

It's important to note that in existing CLTs for LSSs, there often exists an implicit condition stating that the spectrum of $\boldsymbol{\Sigma_n}$ should not degenerate asymptotically. This condition is crucial for the validity and applicability of the CLT, as the degeneracy of the population covariance matrix's spectrum can significantly impact the convergence properties and statistical behavior of the sample covariance matrix. 
The proposed normalization method emerges as a promising approach to address these challenges. By effectively pushing the unbounded spectrum of the population covariance matrix back into a constant order region, our approach simplifies the study of such scenarios. This concept aligns closely with a previous strategy where $\left(\sqrt{p}+\sqrt{n}\right)^2$ was employed to counteract the unbounded spectrum arising from large values of $p$. This paves the way for a more comprehensive understanding and analysis of high-dimensional sample covariance matrix.

The remainder of this paper is organized as follows. In the upcoming section, we will articulate our principal theoretical findings pertaining to the LSD and the CLT for LSS of the renormalized sample covariance matrix. Following that, we will provide an illustrative example demonstrating the application of our results in the context of high-dimensional statistical inference. The proofs of the main theorems will be deferred to Section \ref{proofofmt} before lemmas and a concluding discussion. For the sake of conciseness, some specific computational details will be relegated to the Appendix.

\section{Main results}
This section is dedicated to elucidating our key theoretical findings concerning the spectral properties of the renormalized sample covariance matrix, extending our analysis to encompass varying rates of dimension and sample size without restrictions. Our primary emphasis will be on the convergence of the ESD. Following this, we will introduce a harmonic version of the CLT for LSS specifically designed for the former mentioned scenarios, ensuring its adaptability to a diverse range of cases relevant to statistical applications. 

It is crucial to underscore that our findings naturally hold an intrinsic connection with those established under the original well-studied model. This connection arises from the fact that our renormalization process inherently involves a rescaling of the spectrum of the original sample covariance matrix. Moreover, this strategic renormalization ensures the applicability of commonly used strategies across all the cases we consider.
Nevertheless, it is also important to emphasize that these inherent connections and enhanced applicability do not diminish the necessity and importance of a renewed proof. The influence of the renormalization process warrants a careful and thorough investigation to prevent the inadvertent neglect of contribution terms.

\subsection{Convergence of the ESD}
Let us initiate our exploration by examining the convergence of the ESD of the renormalized sample covariance matrix $$\mathbf S_n=\frac{1}{\left(\sqrt{p}+\sqrt{n}\right)^2\|\boldsymbol{\Sigma}_n\|}{\mathbf X_n}{\mathbf X_n^*}.$$ In this analysis, we extend our consideration beyond the ICS model to a more general framework. Specifically, we allow for nonlinear dependencies among the variables of the population, akin to the scenario explored in \cite{BaiZ08L} under the comparable scenario.
\begin{thm}\label{th:1}
If, as $p\to\infty$ and $n\to\infty$, the following conditions hold:
\begin{enumerate}
	\item [(a)] Assumption on variables:  $\left\{x_{jk}\right\}$ are an array of complex random variables. Let $\mathbf X_n= \left(x_{jk}\right)_{p\times n}$. Assume that the columns $\bx_k,1\le k\le n$ of $\bX_n$ are i.i.d. and $\re(\bx_1\bx_1^*)=\bSig_n$; 
	\item [(b)]Assumption on dependent structure I: 
	 $H_n=F^{\bSig_n/\|\bSig_n\|}$ tends to a non-random probability distribution $H$;
	\item [(c)] Assumption on dependent structure II: For any non-random $p\times p$ matrix $\bB$ with bounded norm and $q\ge2$,
	$$\re\left|\bx_1^*\bB\bx_1-\rtr\(\bB\bSig_n\)\right|^q=O\(p^{q/2}\);$$
	\item [(d)] Asymptotic framework: $c_{n1}=p/\nu_{np}\to c_1 \in [0,1]$ and $c_{n2}=n/\nu_{np}\to c_2\in [0,1]$, where $\nu_{np}=\(\sqrt p+\sqrt n\)^2$.
\end{enumerate}

Then, with probability one, the ESD of the renormalized sample covariance matrix $\mathbf S_n$ converges to a probability distribution $F$, whose Stieltjes transform $m(z)$ satisfies
\begin{align}\label{equa}
	&m(z)=\int\frac1{\(c_2-c_1-c_1zm(z)\)t-z}dH(t)
\end{align}
and is unique in the set $\left\{m\in\mathbb C^+:-(c_2-c_1)/z+c_1m\in\mathbb C^+\ {\rm or}\ c_2=0\right\}$.
\end{thm}

\begin{remark}\label{rem1}
If $\bX_n=\bGa_n\bY_n$, and the entries $y_{jk}$ of $\bY_n$ are i.i.d. such that $\re(y_{jk})=0$, $\re(y_{jk}^2)=1$, $\bGa_n\bGa_n^*=\bSig_n$ , then the condition $(c)$ of Theorem \ref{th:1} holds. We present a verification of this remark in Appendix \ref{verirem2}.
\end{remark}

It is imperative to acknowledge that the LSD is contingent on the values of $c_1$, $c_2$, and $H$. To signify this dependence, we shall use $F^{c_1,c_2,H}$ to represent $F$. Accordingly, the finite sample version $F^{c_{n1},c_{n2},H_{n}}$ is derived from the equation above by substituting ${c_1,c_2,H}$ with $c_{n1},c_{n2},H_{n}$. This substitution is considered safe as we meticulously retain all terms related to the ratio of powers of $p$ and $n$, while selectively neglecting negligible error terms.

The following observations naturally follow from the above theorem, providing a coherent alignment with our anticipated results:
 \begin{itemize}
 	\item In a comparable scenario where $c_1,c_2\in (0,1)$,  the  LSD $F$ is a transformation of the classical M-P law.
 	\item In large $n$ scenario, we have $c_1=0, c_2=1$ thus  $m(z)=\int\frac1{t-z}dH(t)$. In this case $F{=}H.$
 	\item In large $p$ scenario, we have $c_1=1, c_2=0$ thus $m(z)=\int\frac1{(-1-zm(z))t-z}dH(t)$.  This leads to $m(z)=\frac{1}{-z}$  and consequently, $F(x){=}I(x\geq 0).$
 \end{itemize}

\subsection{A harmonic CLT for LSS}
Now, let's delve into the CLT for LSS within the classical ICS model, incorporating a Lindeberg-type condition. For this purpose, we express $\underline{m}(z)=-(c_2-c_1)/z+c_1m(z)$ and introduce a first-stage centralized random process with a rescale, defined as follows:
\begin{align*}
	G_n(x)=\sqrt{\frac p n}{\nu_{np}}\(F^{\bS_n}(x)-F^{c_{n1},c_{n2},H_n}(x)\).
\end{align*}
 We establish the following assumptions:
\begin{itemize}
	\item [(a')] Assumption on variables: For each $n$, let $y_{jk}=y_{jk}^{(n)}$, where $j\le p$ and $k\le n$. These variables are independent, satisfying $\re(y_{jk})=0$, $\re|y_{jk}|^2=1$, and $\re y_{jk}^2=\alpha_y$ where $\alpha_y$ is given by: $\alpha_y=\begin{cases}
		1,\quad {\rm if \ } y_{jk}\ {\rm is\  real}\\
	0,	\quad {\rm if \ } y_{jk}\ {\rm is\  complex}\end{cases}.$ The fourth moment $\re|y_{jk}|^4=2+\alpha_y+\Delta_y.$ Additionally, for any fixed $\delta>0$, it is required that
$
		\frac1{pn}\sum_{jk}\re|y_{jk}|^4I\(|y_{jk}|\ge\delta\sqrt[4]{pn}\)\to0.
$

	\item [(e)] ICS model: $\bX_n=\bGa_n\bY_n=\(\bx_1,\cdots,\bx_n\)$ with $\bSig_n=\bGa_n\bGa_n^*$;
	\item [(f)] Non-trivial support: The support set of $H$, denoted as $\mathcal{S}_{H}$, satisfies that $\mathcal{S}_{H} \cap \mathbb{R}^+$ is a non-empty set.

\end{itemize}
Let's also define several quantities as follows: \begin{align*}
	&g(z)=c_1\int \frac {t}{1+t\underline m(z)}dH(t)-z,\quad h(z_1,z_2)=\int\frac {t^2dH(t)}{\(1+t\underline m(z_1)\)\(1+t\underline m(z_2)\)}, \\&t(z_1,z_2)=\lim_{p}\frac{1}{p}\rtr\left({\widetilde\bGa_n^*}\(\underline m(z_1){\widetilde\bSig_n}+\bI_p\)^{-2}\widetilde\bGa_n\right)\circ\left(\widetilde\bGa_n^*\(\underline m(z_2){\widetilde\bSig_n}+\bI_p\)^{-2}\widetilde\bGa_n\right),\\ & s(z)=\lim_{p}\frac{1}{p}\rtr\(\widetilde\bGa_n^*\(\underline m(z){\widetilde\bSig_n}+\bI_p\)^{-1}\widetilde\bGa_n\)\circ\({\widetilde\bGa_n^*}\(\underline m(z){\widetilde\bSig_n}+\bI_p\)^{-2}\widetilde\bGa_n\),
\end{align*}
where $\widetilde\bGa_n=\frac{\bGa_n}{\sqrt{\|\bSig_n\|}}, \widetilde\bSig_n=\frac{\bSig_n}{\|\bSig_n\|}.$
Then consider $M$ test functions $f_1, \ldots, f_{M}$ that are analytic on an open region containing the interval $[0,1]$. Define $M$ dimensional random vector	\begin{align}\label{cltig}
	\(X_{f_{1n}},\cdots,X_{f_{Mn}}\)=p\(\int f_1(x)dG_n(x),\cdots,\int f_MdG_n(x)\).
	\end{align} With these definitions in place, we aim to establish the following theorem.

\begin{thm}\label{cltthm}
	Assuming the conditions (a'), (b), (d)-(f), the centralized random vector $$\(X_{f_{1n}}-\re X_{f_{1n}},\cdots,X_{f_{Mn}}-\re X_{f_{Mn}}\),$$ with  \begin{align}\label{meanf}
	\re\(X_{f_{jn}}\)=-\frac {1}{2\pi i}\oint_{\mathcal{C}} {f_j(z)}\mathcal{A}_n(z)dz, \quad 1\leq j\leq M,
\end{align} converges weakly to a zero-mean Gaussian vector $\(X_{f_1},\cdots,X_{f_M}\)$ with a covariance matrix $\mathcal{O}=(o_{jk})_{M\times M}$ given by \begin{align}\label{covf}
	o_{jk}={\rm Cov}\(X_{f_j},X_{f_k}\)=-\frac {1}{4\pi^2 }\oint_{\mathcal{C}_1}\oint_{\mathcal{C}_2} {f_j(z)}{f_k(z)}\mathcal{B}\big(z_1,z_2\big)dz_2dz_1.
\end{align}
Here, the contours in \eqref{meanf} and \eqref{covf} are closed and taken in the positive direction in the complex plane, each enclosing the support of $F^{c_1, c_2, H}$ and the parameters
	\begin{align*}
		 \mathcal{A}_n(z)=
\frac{\alpha_y\sqrt {c_{n1}c_{n2}}\(\underline m_n^0(z)\)^3\int\frac{t^2}{\(1+t\underline m(z)\)^3}dH(t)}{c_{n2}^3\(1-\frac{c_{n1}}{c_{n2}}\int\frac{\(\underline m_n^0(z)t\)^2}{\(1+t\underline m_n^0(z)\)^2}dH(t)\)^2}+\frac{\Delta_y\sqrt {c_{n1}c_{n2}}\(\underline m_n^0(z)\)^3s(z) }{c_{n2}^3\(1-\frac{c_{n1}}{c_{n2}}\int\frac{\(\underline m_n^0(z)t\)^2}{\(1+t\underline m_n^0(z)\)^2}dH(t)\)}
	\end{align*}
and 
\begin{align*}
	\mathcal{B}\big(z_1,z_2\big)=&\frac{1+\alpha_y}{\(z_1-z_2\)^2}\Bigg[\frac{c_1c_2\(h^2(z_1,z_2)-h(z_1,z_1)h(z_2,z_2)\)}{\left[g^2(z_1)-c_1c_2h(z_1,z_1)\right]\left[g^2(z_2)-c_1c_2h(z_2,z_2)\right]}\\
	&+\frac{g^2(z_1)h(z_2,z_2)+h(z_1,z_1)g^2(z_2)-2g(z_1)g(z_2)h(z_1,z_2)}{\left[g^2(z_1)-c_1c_2h(z_1,z_1)\right]\left[g^2(z_2)-c_1c_2h(z_2,z_2)\right]}\Bigg]\\
	&+\frac{\Delta_yt(z_1,z_2)}{\(g^2(z_1)-c_1c_2h(z_1,z_1)\)\(g^2(z_2)-c_1c_2h(z_2,z_2)\)}.
\end{align*}
 
\end{thm}

{ In the following discussion, we provide some crucial insights.} It is evident from the expression in the preceding theorem that the mean and variance-covariance of the LSSs converge in all cases under our assumptions, irrespective of the divergent rates of $p$ and $n$, as well as the order of the spectral norm of the population covariance matrix. The distinction lies in the fact that the limit of the mean function remains of constant order in the comparable scenario but degenerates to 0 under the other two scenarios. Conversely, the limit of the variance-covariance function consistently remains of constant order. Importantly, in practical applications, especially under non-Gaussian scenarios, practitioners often resort to finite sample versions of these statistics. Consequently, this phenomenon does not hinder the utility of our CLT, as indicated in the section of application.

On the flip side, we aim to glean insights into the convergence rate of specific LSS associated with the original sample covariance matrix. For a deeper understanding, consider the normalization parameter for $G_n(x)$. From this, one can deduce that the convergence rate of  $\rtr\bS_n^0$ follows the order of  $\frac{\sqrt{np}}{\nu_{n,p}\|\boldsymbol{\bSig}_n\|}$ whereas the convergence rate of $\rtr\(\bS_n^0\)^2$ is of the order $\frac{n\sqrt{np}}{\nu_{n,p}^2\|\boldsymbol{\bSig}_n\|^2}$. An overarching insight emerges, highlighting that the convergence rate of LSS is significantly influenced by the rates of $p$ and $n$, as well as the spectral norm $\|\boldsymbol{\Sigma}_n\|$ and the test function $f(\cdot)$. Notably, these findings align with the results presented in \cite{Yin2022} within a comparable scenario. A more detailed examination of these observations is available in the subsequent application section.

\section{Statistical application: Identity covariance matrix test}

In this section, we showcase the substantial advantages of our established harmonic CLT in enhancing statistical inference, with a particular focus on the classical identity test for the population covariance matrix.

 Let us consider a scenario where samples are drawn from a $p$-dimensional population. Consider a complex data matrix $\mathbf{X}_n = (\mathbf{x}_1, \ldots, \mathbf{x}_n)$ with dimensions $p \times n$, where $n$ i.i.d. $p$-dimensional random vectors $\{\mathbf{x}_j = \boldsymbol{\Sigma}_n^{1/2} \mathbf{y}_j\}_{1 \leqslant j \leqslant n}$. Understanding the covariance structure among these variables is pivotal for developing robust statistical tools, making the test for covariance structure a foundational step in various statistical inferences. Specifically, when the population adheres to the ICS model, the test for a given covariance structure involves transforming the samples and executing the identity test:
\[ H_0: \boldsymbol{\Sigma}_n = \boldsymbol{I}_p \quad \text{vs.} \quad \boldsymbol{\Sigma}_n \neq \boldsymbol{I}_p. \]
This identity test has garnered considerable attention over the years. 
\subsection{Frobenius norm test} In the classical large sample framework, where $p$ is considered fixed while the sample size $n$ tends to infinity, \cite{Nagao73A} proposed utilizing the Frobenius norm test statistics
\[ V = \frac{1}{p} \operatorname{tr}\left\{(\mathbf{S}_n^0 - \mathbf{I}_p)^2\right\}. \]
However, the inadequacy of this test in high-dimensional scenarios, where $p$ is comparable to $n$, was highlighted by \cite{LedoitW02S}. They introduced a modification to Nagao's test statistics as
\[ W = \frac{1}{p} \operatorname{tr}\left\{(\mathbf{S}_n^0 - \mathbf{I}_p)^2\right\} - \frac{1}{n p}\left\{\operatorname{tr}(\mathbf{S}_n^0)\right\}^2 + \frac{p}{n}. \]
Demonstrating, under Gaussian assumptions and in the comparable scenario, that the null distribution of their test statistics is asymptotically Gaussian with
\[ n W - p - 1 \stackrel{d}{\longrightarrow} {N}(0,4). \]
This insight was recently extended by \cite{Qiu2023}, accommodating non-Gaussian underlying distributions with a finite $6+\varepsilon$ order moment. They relaxed the comparable scenario to an ultra-high dimensional case where $ p \wedge n \rightarrow \infty \text { and } n^2 / p=O(1)$. Under their framework, by investigation of the matrix $\frac{1}{\sqrt{n p}}\left(\mathbf{X}_n^{T} \mathbf{X}_n-p \mathbf{I}_n\right)$, the null distribution of $W$ is demonstrated to be
\[ nW - p - (\re y_{11}^4  - 2) \stackrel{d}{\longrightarrow} {N}(0,4), \]
essentially aligning with the classical comparable scenario.

In this section, we harness the power of our established harmonic CLT to showcase that the null distribution of the statistics $W$ is indeed independent of the rates of $p$ and $n$ under a relaxed moment assumption, incorporating the optimal fourth-moment assumption. This in-depth analysis serves to fortify the robustness of Ledoit and Wolf's identity test statistic under ICS. 
To be specific, we establish the following theorem under both real and complex variable cases.  
\begin{thm}\label{identitytestf}
 Assume the Assumption (a').
Under the null hypothesis $H_0$, as $p \wedge n \rightarrow \infty$, the asymptotic distribution of the test statistic $W$ is given by
\begin{align*}
    n W - p - (\alpha_y+\Delta_y) \stackrel{d}{\longrightarrow} N\(0,2\(1+\alpha_y\)\).
\end{align*}
\end{thm}

\begin{remark}
We wish to apprise the reader that the non-null distribution of this statistic can be derived utilizing the same methodology as in the proof of the preceding theorem. Employing a parallel approach allows one to readily establish an asymptotic distribution, as demonstrated in Theorem 4.2 of \cite{Qiu2023}, without constraints on the rate of $p$. This, in turn, enables a similar power analysis to be conducted.
\end{remark}

\subsection{Likelihood ratio test} Another essential statistic in the context of identity tests is the Likelihood Ratio Test (LRT), defined as
\begin{align*}
L^* = \operatorname{tr} \mathbf{S}_n^0 - \log |\mathbf{S}_n^0| - p.
\end{align*}
In the classical scenario where $p$ is fixed and $n \rightarrow \infty$, the well-established theory in multivariate analysis asserts that $nL^*$ converges to the chi-square distribution with degrees of freedom $p(p+1)/2$ under the null hypothesis $H_0$.

However, the limitations of LRT in high-dimensional cases have been exposed by  \cite{BaiJ09C}. To address this, they proposed a correction to LRT based on the CLT for the LSS established by Bai and Silverstein in \cite{BaiS04C}, assuming a comparable scenario with a Gaussian-like fourth-moment assumption. Importantly, their analysis assumes that the limit of the dimension-to-sample-size ratio $p/n$ is less than 1.

Under these assumptions, it is proven that under $H_0$,
\begin{align*}
\frac{L^* - p A_n(c_n) - B_n(c_n)}{C_n(c_n)} \stackrel{d}{\longrightarrow} \mathcal{N}(0,1),
\end{align*}
where $c_n = p/n$ and $A_n(c_n) = 1 - \frac{c_n-1}{c_n} \log (1-c_n)$, $B_n(c_n) = -\frac{\log (1-c_n)}{2}$, and $C^2_n(c_n) = -2 \log (1-c_n) - 2 c_n$. This result is extended in \cite{JiangY13C} to cases where the limit of $p/n$ can be equal to 1 under Gaussian assumptions.

In the work of \cite{Qiu2023}, the extension of LRT to the ultra-high-dimensional case is considered under a Gaussian-like fourth-moment assumption. However, when $p>n$, the LRT statistic is degenerate and undefined due to $\left|\mathbf{S}_n\right|=0$. Consequently, the quasi-LRT test statistic is introduced:
\begin{align*}
L = \operatorname{tr}\widehat{\mathbf{S}}_n^0 - \log \left|\widehat{\mathbf{S}}_n^0\right| - n,
\end{align*}
where $\widehat{\mathbf{S}}_n = \frac{1}{p} \mathbf{X}_n^{*} \mathbf{X}_n$ is obtained by interchanging the positions of $p$ and $n$, rendering it invertible. In the limit where $p \wedge n \rightarrow \infty$ and $n^2 / p = O(1)$, the null distribution of ${L}$ converges to
\begin{align*}
\frac{L - n A_n(1/c_n) - B_n(1/c_n)}{C_n(1/c_n)} \stackrel{d}{\longrightarrow} \mathcal{N}(0,1),
\end{align*}
aligning with the comparable scenario. 

In this section, we aim to remove the restrictive conditions on both the rates of $p$ and $n$ as well as the Gaussian-like fourth-moment assumption by leveraging our established harmonic CLT. Specifically, we present the following theorem.
 \begin{thm}\label{identitytestl}
 Assume the Assumption (a').
Under the null hypothesis $H_0$, as $p \wedge n \rightarrow \infty$, the asymptotic distribution of the test statistic $L$ and $L^*$ is given by
\begin{align*}
\begin{cases}
	\frac{L^* - p A_n(c_n) - B_n(c_n)}{C_n(c_n)} \stackrel{d}{\longrightarrow} {N}(0,1), & p<n, \\
	\frac{L - n A_n(c_n^{-1}) - B_n(c_n^{-1})}{C_n(c_n^{-1})} \stackrel{d}{\longrightarrow} {N}(0,1), & p>n,
\end{cases} 
\end{align*}
where $c_n = p/n$ and $A_n(c_n) = 1 - \frac{c_n-1}{c_n} \log (1-c_n)$, $B_n(c_n) = -\frac{{\alpha_y}\log (1-c_n)}{2}+\frac{\Delta_yc_{n}}{2}$, and $C^2_n(c_n) = -\(1+\alpha_y\) \log (1-c_n) - \(1+\alpha_y\) c_n$. 
\end{thm}

\begin{remark}
The results presented in the theorem above surprisingly unveil a phenomenon. Specifically, the null distribution of the LRT statistic becomes independent of the certain underlying distribution as both $p$ and $n$ become large and exhibit distinct divergent rates. However, it does rely on the fourth moment of the underlying distribution under the comparable scenario.
\end{remark}

\section{Proof of main theorems}\label{proofofmt}
In this section, we embark on the rigorous proof of our theoretical results. 
While our results establish a bridge to well-established models and broaden the scope of applicability for common strategies, a renewed proof remains a critical step. This is not just a formality but a means to rigorously validate and fortify the robustness of our theoretical framework. The impact of renormalization introduces nuances that demand meticulous attention, and a renewed proof becomes indispensable in ensuring the validity and reliability of our theories across diverse scenarios. 

Before we commence, let's revisit the normalization of the sample covariance matrix using $\|\bSig_n\|$. This is equivalent to transforming $\bGa_n$ into $\widetilde\bGa_n$. However, for the sake of simplicity in notation, we would like to emphasize that in the subsequent proofs, we will interchangeably use $\bGa_n$ to denote $\widetilde\bGa_n$, treating $\bGa_n$ as a matrix with unit norm and ESD $H_n$.

Before delving into the proof, we introduce several notations to streamline the subsequent analysis: $\bA(z)=\bS_n-z\bI_p,$
$\bA_k(z)= \bA(z)-\nu_{np}^{-1}\bx _k\bx_k^*,$ $\bA_{kj}(z)= \bA_k(z)-\nu_{np}^{-1}\bx _j\bx_j^*,$
\begin{align*}
	\varepsilon_k(z)=\bx_k^*\bA_k^{-1}(z)\bx_k-\rtr\(\bA_k^{-1}(z)\bSig_n\), \quad
	\eta_k(z)=\bx_k^*\bA_k^{-2}(z)\bx_k-\rtr\(\bA_k^{-2}(z)\bSig_n\),
\end{align*}and
\begin{align*}
\beta_k(z)=\frac1{1+\nu_{np}^{-1}\bx_k^*\bA_k^{-1}(z)\bx_k},\quad b_{nk}(z)=\frac1{1+\nu_{np}^{-1}\rtr\(\bA_k^{-1}(z)\bSig_n\)},\\
\beta_{kj}(z)=\frac1{1+\nu_{np}^{-1}\bx_j^*\bA_{kj}^{-1}(z)\bx_j},\quad b_n(z)=\frac1{1+\nu_{np}^{-1}\re\rtr\(\bA^{-1}(z)\bSig_n\)}.
\end{align*}
It is evident that
$
	\|\bA_k^{-1}(z)\|\le v^{-1},\quad |\beta_k(z)|\le |z|/v,\quad |b_n(z)|\le |z|/v.
$

\subsection{Proof of Theorem \ref{th:1}}

We delineate our proof into three distinct steps, each contributing to the establishment of Theorem \ref{th:1}:

\begin{itemize}
    \item \textbf{Step 1:} For any fixed $z \in \mathbb{C}^+$, we demonstrate the almost sure convergence $m_n(z) - \re(m_n(z)) \xrightarrow{\text{a.s.}} 0$.
    
    \item \textbf{Step 2:} Further, for any fixed $z \in \mathbb{C}^+$, we establish the convergence $\re(m_n(z)) - m(z) \rightarrow 0$.
    Leveraging Lemma \ref{lemma:12}, we conclude that for every $z \in \mathbb{C}^+$,
    \begin{align*}
        m_n(z) - m(z) \xrightarrow{\text{a.s.}} 0.
    \end{align*}
    
    \item \textbf{Step 3:} For $z \in \mathbb{C}^+$, we assert that $m(z) \in \mathbb{C}^+$ is the unique solution to the equation \eqref{equa}.
\end{itemize}
After completing these three steps, the result of Theorem \ref{th:1} follows.
The details of the proof for these steps will be deferred to the Appendix \ref{Dtlsd}.

\subsection{Proof of Theorem \ref{cltthm}}
To establish the CLT, we adopt the strategy employed in \cite{BaiS10S}, which combines the Stieltjes transform method with the decomposition of martingale difference sequences, followed by the application of the CLT for martingales. We also direct the reader to \cite{ZhengB15S, Hu2019, Qiu2023, GaoH17H}, and references therein, for instances where a similar strategy has been applied. Consider the Cauchy Integral Formula:
\begin{align*}
	\int f(x) \, dG(x) = -\frac{1}{2\pi i} \oint_{\mathcal{C}_G} f(z) m_G(z) dz.
\end{align*}
This formula is valid for any distribution function $G$ and any $f$ analytic on an open set containing the support of $G$. The complex integral on the right is over any positively oriented contour enclosing the support of $G$ and on which $f$ is analytic. Thus, we shift our focus from investigating the LSS to examining the convergence of the normalized Stieltjes transform of the ESD of our renormalized sample covariance matrix $\bS_n$.
Let
\begin{align*}
	M_n(z)=\sqrt{\frac pn}\nu_{np}\left[m_n(z)-m_n^0(z)\right]=\sqrt{\frac pn}\nu_{np}\left[m_{F^{\bS_n}}(z)-m_{F^{c_{n1},c_{n2},H_n}}(z)\right],
\end{align*}
and
\begin{align*}
	&m_n^0(z)=-z^{-1}\int\frac1{\underline m_n^0(z)t+1}dH_n(t),\quad \underline m_n^0(z)=-\(c_{2n}-c_{1n}\)/z+c_{1n}m_n^0(z).
\end{align*}
Then the outline of the proof is organized as follows.

\begin{description}
	\item \textbf{Step 1:} We commence with a truncation step to limit the underlying variables to a suitable order, depending on $p$ and $n$, thereby avoiding the need to handle excessively high-order moments.
	\item \textbf{Step 2:} Let $v_0 > 0$ be arbitrary. Consider $x_r$ as any number greater than 1 and $x_{\ell}$ as any negative number. Define $\mathcal C_u=\{x+iv_0:x\in[x_{\ell},x_r]\}$ and $\mathcal C_{\ell}=\{x_{\ell}+iv:v\in[0,v_0]\}$. Set $\mathcal C=\mathcal C^+\cup\overline{\mathcal C^+}$, where $\mathcal C^+=C_{\ell}\cup\mathcal C_u\cup\{x_{r}+iv:v\in[0,v_0]\}$.
Moreover, define the subsets $\mathcal C_n$ of $\mathcal C^+$. Choose a sequence $\{\varepsilon_n\}$ decreasing to zero, satisfying, for some $\alpha\in(0,1)$, $\varepsilon_n\ge n^{-\alpha}$. Let $\mathcal C_r=\{x_{r}+iv:v\in[n^{-1}\varepsilon_n,v_0]\}$; then, $\mathcal C_n=C_{\ell}\cup\mathcal C_u\cup\mathcal C_r$.
Under the conditions of Theorem \ref{cltthm}, we ascertain that, for any $f\in\{f_1,\cdots,f_M\}$ and for all large $n$, with probability 1,
\begin{align*}
	\int f(x)dG_n(x)=-\frac1{2\pi i}\oint_{\mathcal C}f(z)M_n^0(z)dz.
\end{align*} 
To eliminate the small probability event where some sample eigenvalues fall outside the defined interval, we introduce a truncated version of the process $M_n^0$. The truncated process is defined as follows, for $z=x+iv$:
\begin{align*}
	 M_n(z)=\begin{cases}
		M_n^0(z),\qquad\qquad\qquad \text{for}\ z\in\mathcal C_n,\\
		M_n^0(x_r+in^{-1}\varepsilon_n),\quad\ \ \text{for}\ x=x_r, v\in[0,n^{-1}\varepsilon_n].
	\end{cases}
\end{align*}
It can be verified that, with probability 1,
\begin{align*}
	\oint_{\mathcal C}f(z)\(M_n(z)-M_n^0(z)\)dz\to0.
\end{align*}
Subsequently, we establish the CLT for the process $M_n(z)$ by decompose it into two parts as \begin{align*}
	M_{n1}(z)=\sqrt{\frac pn}\nu_{np}[m_n(z)-\re m_n(z)],\quad
	M_{n2}(z)=\sqrt{\frac pn}\nu_{np}[\re m_n(z)- m_n^0(z)].
\end{align*}
\item \textbf{Step 3:} Prove the tightness of $M_{n1}(z)$ and finish the proof.
\end{description}

We will now systematically advance through each of these steps.

\subsubsection{Truncation and recentralization}
By condition $(a')$, we can choose $\delta_n$ decreasing to 0 such that
$
	\frac1{np\delta_n^4}\sum_{jk}\re|y_{jk}|^4I\(|y_{jk}|\ge\delta_n\sqrt[4]{pn}\)\to0.
$
Let
$
	\mathbb S_n=\frac1{\nu_{np}}\bGa_n\mathbb Y_n\mathbb Y_n^*\bGa_n^* 
$
with $\mathbb Y_n$ being a $p\times n$ matrix having ($j,k$)-th entry $y_{jk}I\(|y_{jk}|<\delta_n\sqrt[4]{pn}\)$. Then it yields
\begin{align*}
	{\rm P}\(\bS_n\neq\mathbb S_n\)\le\sum_{jk}{\rm P}\(|y_{jk}|\ge \delta_n\sqrt[4]{pn}\)
	\le\frac1{pn\delta_n^4}\sum_{jk} \re\left[|y_{jk}|^4I\(|y_{jk}|\ge \delta_n\sqrt[4]{pn}\)\right]\to0.
\end{align*}
Define $\hat{\mathbb S}_n=\frac1{\nu_{np}}\bGa_n\hat{\mathbb Y}_n\hat{\mathbb Y}_n^*\bGa_n^*$ with $\hat{\mathbb Y}_n$ being a $p\times n$ matrix having ($j,k$)-th entry $$\left[y_{jk}I\(|y_{jk}|<\delta_n\sqrt[4]{pn}\)-\re y_{jk}I\(|y_{jk}|<\delta_n\sqrt[4]{pn}\)\right]/\sigma_{jk},$$
where $\sigma_{jk}^2={\rm Var}\(y_{jk}I\(|y_{jk}|<\delta_n\sqrt[4]{pn}\)\)$. We use $\mathbb G_n(x)$ and $\hat{\mathbb G}_n(x)$ to denote the analogues of $G_n(x)$ with the matrix $\bS_n$ replaced by $\mathbb S_n$ and $\hat{\mathbb S}_n$, respectively. Using the mean value theorem, we have
\begin{align*}
	&\re\left|\int f_j(x)d\mathbb G_n(x)-\int f_j(x)d\hat{\mathbb G}_n(x)\right|\le \frac{C\nu_{np}}{\sqrt{pn}}\sum_{k=1}^p\re\left|\lambda_k^{\mathbb S_n}-\lambda_k^{\hat{\mathbb S}_n}\right|\\
\le&\frac C{\sqrt{pn}} \(\re\rtr\(\mathbb Y_n-\hat{\mathbb Y}_n\)\(\mathbb Y_n-\hat{\mathbb Y}_n\)^*\)^{1/2}\(\re\rtr\(\mathbb Y_n\mathbb Y_n^*+\hat{\mathbb Y}_n\hat{\mathbb Y}_n^*\)\)^{1/2}\\
\le&C\(\sum_{jk}\left|\re y_{jk}I\(|y_{jk}|<\delta_n\sqrt[4]{pn}\)\right|^2+\sum_{jk}\(1-\sigma_{jk}^2\)^2\)^{1/2}\to 0,
\end{align*}
where $\lambda_k^{\bA}$ denotes the $k$-th smallest eigenvalues of $\bA$.
Therefore, we obtain
\begin{align*}
	\int f_j(x)dG_n(x)=\int f_j(x)d\hat{\mathbb G}_n(x)+o_p(1).
\end{align*}
For simplicity, in the following, the truncated and centralized variables are still denoted by $y_{jk}$, and we assume that
\begin{align*}
	|y_{jk}|\le \delta_n\sqrt[4]{pn},\ \re\(y_{jk}\)=0,\ \re|y_{jk}|^2=1,\ \re|y_{jk}|^4=\Delta_y+2+\alpha_y,
\end{align*}
and
$\re y_{jk}^2=\alpha_y=\begin{cases}
	1,\quad y_{jk}\ {\rm is\  real},\\
	0,	 \ \ \ \  y_{jk}\ {\rm is\  complex}.\end{cases}$

\subsubsection{Convergence of Stieltjes transform}
Under the assumptions of Theorem \ref{cltthm}, and after the truncation, we will establish the following lemma, which informs us about the asymptotic properties of the random process $M_n(z)$.\begin{lemma}\label{CLTST}
	Under conditions of Theorem \ref{cltthm}, $\{ M_n(z)-\mathcal{A}_n(z)\}$ forms a tight sequence on $\mathcal C^+$ and converges weakly to a two-dimensional zero mean Gaussion process $M(\cdot)$ with a covariance function
\begin{align*}
	{\rm Cov}\big(M&(z_1),M(z_2)\big)=\frac{1+\alpha_y}{\(z_1-z_2\)^2}\Bigg[\frac{c_1c_2\(h^2(z_1,z_2)-h(z_1,z_1)h(z_2,z_2)\)}{\left[g^2(z_1)-c_1c_2h(z_1,z_1)\right]\left[g^2(z_2)-c_1c_2h(z_2,z_2)\right]}\\
	&+\frac{g^2(z_1)h(z_2,z_2)+h(z_1,z_1)g^2(z_2)-2g(z_1)g(z_2)h(z_1,z_2)}{\left[g^2(z_1)-c_1c_2h(z_1,z_1)\right]\left[g^2(z_2)-c_1c_2h(z_2,z_2)\right]}\Bigg]\\
	&+\frac{\Delta_yt(z_1,z_2)}{\(g^2(z_1)-c_1c_2h(z_1,z_1)\)\(g^2(z_2)-c_1c_2h(z_2,z_2)\)}.
\end{align*}
Here the mean parameter
	\begin{align*}
		\mathcal{A}_n(z)=
\frac{\alpha_y\sqrt {c_{n1}c_{n2}}\(\underline m_n^0(z)\)^3\int\frac{t^2}{\(1+t\underline m(z)\)^3}dH(t)}{c_{n2}^3\(1-\frac{c_{n1}}{c_{n2}}\int\frac{\(\underline m_n^0(z)t\)^2}{\(1+t\underline m_n^0(z)\)^2}dH(t)\)^2}+\frac{\Delta_y\sqrt {c_{n1}c_{n2}}\(\underline m_n^0(z)\)^3s(z) }{c_{n2}^3\(1-\frac{c_{n1}}{c_{n2}}\int\frac{\(\underline m_n^0(z)t\)^2}{\(1+t\underline m_n^0(z)\)^2}dH(t)\)}.
	\end{align*}

\end{lemma}
To prove this lemma, we write for $z\in\mathcal C_n$, $M_n(z)=M_{n1}(z)+M_{n2}(z)$, where
\begin{align*}
	M_{n1}(z)=\sqrt{\frac pn}\nu_{np}[m_n(z)-\re m_n(z)],\quad
	M_{n2}(z)=\sqrt{\frac pn}\nu_{np}[\re m_n(z)- m_n^0(z)].
\end{align*}
Next, our goal is to determine the limiting distribution of $M_{n1}(z)$ and investigate the asymptotics of $M_{n2}(z)$.  
\paragraph{The limiting distribution of $M_{n1}(z)$}
In this section, our primary objective is to determine the limiting distribution of $M_{n1}(z)$. Specifically, we aim to demonstrate that for any positive integer $r$, the sum
$$\sum_{j=1}^r\alpha_jM_{n1}(z_j) \qquad\Im z_j\neq0$$
converges in distribution to a Gaussian random variable. To achieve this, we will leverage the central limit theorem for martingale difference sequences. Since
\begin{align*}
	\lim_{v_0\downarrow 0}\limsup_{n\to\infty}\re\left|\int_{\mathcal{C}_l\cup\mathcal{C}_r}f(z)M_{n1}(z)dz\right|^2\to0,
\end{align*}
it suffices to consider $z=u+iv_0\in \mathcal{C}_u$. 
By Section \ref{se:1} and 
\begin{align*}
	\beta_k(z)=&b_{nk}(z)-\frac1{\nu_{np}}\beta_k(z)b_{nk}(z)\varepsilon_k(z)\\
	=&b_{nk}(z)-\frac1{\nu_{np}}b_{nk}^2(z)\varepsilon_k(z)+\frac1{\nu_{np}^2}\beta_k(z)b_{nk}^2(z)\varepsilon_k^2(z),
\end{align*}
it is obvious that
\begin{align*}
	M_{n1}(z)=&-\frac{1}{\sqrt {pn}}\sum_{k=1}^n\re_kb_{nk}(z)\eta_k(z)+\frac{1}{\nu_{np}\sqrt {pn}}\sum_{k=1}^n\re_kb_{nk}^2(z)\varepsilon_k(z)\rtr\(\bA_k^{-2}(z)\bSig_n\)\\
	&+\frac{1}{\nu_{np}\sqrt {pn}}\sum_{k=1}^n\(\re_k-\re_{k-1}\)\beta_k(z)b_{nk}(z)\varepsilon_k(z)\eta_k(z)\\
	&-\frac{1}{\nu_{np}^{2}\sqrt {pn}}\sum_{k=1}^n\(\re_k-\re_{k-1}\)\beta_k(z)b_{nk}^2(z)\varepsilon_k^2(z)\rtr\(\bA_k^{-2}(z)\bSig_n\).
\end{align*}
By Lemma \ref{lep1} and Cauchy-Schwarz inequality, we find
\begin{align*}
	&\re\left|\frac{1}{\nu_{np}\sqrt {pn}}\sum_{k=1}^n\(\re_k-\re_{k-1}\)\beta_k(z)b_{nk}(z)\varepsilon_k(z)\eta_k(z)\right|^2\\
\le&\frac{C}{\nu_{np}^2pn}\sum_{k=1}^n\re^{1/2}\left|\varepsilon_k(z)\right|^4\re^{1/2}\left|\eta_k(z)\right|^4\le C\delta_n^4\to0,
\end{align*}
and
\begin{align*}
	&\re\left|\frac{1}{\nu_{np}^{2}\sqrt {pn}}\sum_{k=1}^n\(\re_k-\re_{k-1}\)\beta_k(z)b_{nk}^2(z)\varepsilon_k^2(z)\rtr\(\bA_k^{-2}(z)\bSig_n\)\right|^2\\
	\le&\frac{C}{\nu_{np}^3n}\sum_{k=1}^n\re\left|\varepsilon_k(z)\right|^4\le C\delta_n^4\to0.
\end{align*}
Consequently, we obtain
\begin{align*}
	M_{n1}(z)=&-\frac{1}{\sqrt {pn}}\sum_{k=1}^n\re_k\frac{d}{dz}b_{nk}(z)\varepsilon_k(z)+o_p(1)\triangleq \sum_{k=1}^nY_k(z)+o_p(1).
\end{align*}
Therefore, it is sufficient to demonstrate $$\sum_{j=1}^r\alpha_j\sum_{k=1}^nY_k(z_j)=\sum_{k=1}^n\sum_{j=1}^r\alpha_jY_k(z_j)$$ converges in distribution to a Gaussian random variable. By applying Lemma \ref{lep2}, we need to verify conditions $(i)$ and $(ii)$. Utilizing Lemma \ref{lep1}, we obtain \ref{lep1}, we have
\begin{align*}
	\sum_{k=1}^n\re\left|\sum_{j=1}^r\alpha_jY_k(z_j)\right|^4\le&\frac C{\nu_{np}^4}\sum_{k=1}^n\sum_{j=1}^r\alpha_j^4\bigg[\re|\gamma_k(z_j)|^4+
	\re|\varepsilon_k(z_j)|^4\bigg]\le  C\delta_n^4\to0,
\end{align*}
which implies that condition $(ii)$ of Lemma \ref{lep2} is satisfied. Subsequently, we need to find the limit in probability of $$\Phi(z_1,z_2)\triangleq\sum_{k=1}^n\re_{k-1}\left[Y_k(z_1)Y_k(z_2)\right]$$
for $z_1,z_2$ with nonzero imaginary parts. It is obvious that
\begin{align}\label{inp}
	\Phi(z_1,z_2)=\frac1{pn}\frac{\partial^2}{\partial z_2\partial z_1}\sum_{k=1}^n\re_{k-1}\left[\re_k\left(b_{nk}(z_1)\varepsilon_k(z_1)\right)\re_k
	\left(b_{nk}(z_2)\varepsilon_k(z_2)\right)\right].
\end{align}
Looking at the analysis on page 571 in \cite{BaiS10S}, it is sufficient to prove that the following condition holds:
\begin{align}\label{inp1}\begin{split}
&\frac1{pn}\sum_{k=1}^n\re_{k-1}\left[\re_k\left(b_{nk}(z_1)\varepsilon_k(z_1)\right)\re_k
\left(b_{nk}(z_2)\varepsilon_k(z_2)\right)\right]\\
=&\frac{b_{n}(z_1)b_{n}(z_2)}{pn}\sum_{k=1}^n\re_{k-1}\left[\re_k\left(\varepsilon_k(z_1)\right)\re_k
\left(\varepsilon_k(z_2)\right)\right]+o_p(1)\end{split}\end{align}
converges in probability to a constant.  In this context, we utilized the relationship 
\begin{align}\label{bnk}
	\re|b_{nk}(z)-b_n(z)|^2\le \frac{C}{\nu_{np}}\to0.
\end{align}
Therefore, our goal turns to find the limit in probability of
\begin{align}\label{inp2}
	&\frac{b_{n}(z_1)b_{n}(z_2)}{pn}\sum_{k=1}^n\re_{k-1}\left[\re_k\left(\varepsilon_k(z_1)\right)\re_k
	\left(\varepsilon_k(z_2)\right)\right]\notag\\
=&\frac{b_{n}(z_1)b_{n}(z_2)}{pn}\sum_{k=1}^n\bigg[\(1+\alpha_y\)\rtr\(\re_k\bA_k^{-1}(z_1)\bSig_n\re_k\bA_k^{-1}(z_2)\bSig_n\)\\
&\qquad\qquad\quad+\Delta_y\rtr\left(\bGa_n^*\re_k\bA_k^{-1}(z_1)\bGa_n\right)\circ\left(\bGa_n^*\re_k\bA_k^{-1}(z_2)\bGa_n\right)\bigg].\notag
\end{align}
To this end, we write
$
\bbH_n(z)=\frac {n-1}{\nu_{np}}b_n(z)\bSig_n-z\bI_p.
$
Subsequently, analogous to \eqref{hbound}, we obtain $\left|\bbH_n^{-1}(z)\right|\le C/{v_0^2}.$ It is noteworthy that from \eqref{inver}, we have
\begin{align*}
	\bA_k^{-1}(z)-\bbH_n^{-1}(z)
	=&-\frac{b_n(z)}{\nu_{np}}\sum_{j\neq k}^n\bbH_n^{-1}(z)\(\bx_j\bx_j^*-\bSig_n\)\bA_{kj}^{-1}(z)\\
	&-\frac1{\nu_{np}}\sum_{j\neq k}^n\(\beta_{kj}(z)-b_n(z)\)\bbH_n^{-1}(z)\bx_j\bx_j^*\bA_{kj}^{-1}(z)\\
	&-\frac{b_n(z)}{\nu_{np}^2}\sum_{j\neq k}^n\beta_{kj}(z)\bbH_n^{-1}(z)\bSig_n\bA_{kj}^{-1}(z)\bx_j\bx_j^*\bA_{kj}^{-1}(z)\\
	\triangleq&\bB_1(z)+\bB_2(z)+\bB_3(z).
\end{align*}
Consider ${\bf M}$ as a $p \times p$ matrix with a nonrandom bound on the spectral norm, applicable for all parameters governing ${\bf M}$ and under all realizations of ${\bf M}$. Applying Lemma \ref{lep1}, \eqref{inver}, and \eqref{bnk}, we can conclude that
\begin{align*}
	\re|\beta_{kj}(z)-&b_n(z)|^2\le2\re\left|\beta_{kj}(z)-b_{nk}(z)\right|^2+2\re\left|b_{nk}(z)-b_n(z)\right|^2\\
\le&\frac C{\nu_{np}^2}\re\left|\bx_j^*\bA_{kj}^{-1}(z)\bx_j-\rtr\(\bA_{kj}^{-1}(z)\bSig_n\)\right|^2+\frac C{\nu_{np}}=\frac C{\nu_{np}}.	
\end{align*}
Therefore, we can infer that
\begin{align*}
	\re\bigg|\rtr\(\bB_2(z)+\bB_3(z)\)&{\bf M}\bigg|\le\frac{Cn}{\nu_{np}^{3/2}}\rtr^{1/2}\({\bf M}{\bf M}^*\)+\frac{Cn}{\nu_{np}^{3/2}}\left|\rtr\bbH_n^{-1}(z)\bSig_n\bA_{kj}^{-1}(z){\bf M}\right|\\
	&+\frac{Cn}{\nu_{np}^{2}}\left|\rtr\bbH_n^{-1}(z)\bSig_n\bA_{kj}^{-1}(z)\bSig_n\bA_{kj}^{-1}(z){\bf M}\right|
\end{align*}
and
$
	\re\bigg|\rtr\(\bB_1(z){\bf M}\)\bigg|\le\frac{C\sqrt pn}{\nu_{np}^{3/2}}.
$
Furthermore, for a nonrandom vector $\bu$ with $\bu^*\bu \leq C$, 
\begin{align*}
	\re\bigg|\bu^*\bB_1(z)\bu\bigg|^2=&\frac{b_n^2(z)}{\nu_{np}^2}\sum_{j\neq k}^n\re\bigg|\bu^*\bbH_n^{-1}(z)\(\bx_j\bx_j^*-\bSig_n\)\bA_{kj}^{-1}(z)\bu\bigg|^2\\
	&+\frac{b_n^2(z)}{\nu_{np}^2}\sum_{j_1\neq j_2}^n\re\bigg(\bu^*\bbH_n^{-1}(z)\(\bx_{j_1}\bx_{j_1}^*-\bSig_n\)\bA_{kj_1}^{-1}(z)\bu\\
	&\qquad\cdot\bu^*\bA_{kj_2}^{-1}(\bar z)\(\bx_{j_2}\bx_{j_2}^*-\bSig_n\)\bbH_n^{-1}(\bar z)\bu\bigg)\to0.
\end{align*}
Hence, it follows that
\begin{align}\label{con1}
	&\frac{b_{n}(z_1)b_{n}(z_2)}{p}\left[{z_1}\rtr\(\re_k\bA_k^{-1}(z_1)\bSig_n\)-{z_2}\rtr\(\re_k\bA_k^{-1}(z_2)\bSig_n\)\right]\notag\\
	=&\frac{b_{n}(z_1)b_{n}(z_2)}{p}\left[{z_1}\rtr\(\bbH_n^{-1}(z_1)\bSig_n\)-{z_2}\rtr\(\bbH_n^{-1}(z_2)\bSig_n\)\right]+o_p(1)\notag\\
=&\frac{z_1z_2\(\underline m_n^0(z_1)-\underline m_n^0(z_2)\)}{c_{n1}c_{n2}}\left[1-\frac{\underline m_n^0(z_1)\underline m_n^0(z_2)\(z_1-z_2\)}{c_{n2}\(\underline m_n^0(z_1)-\underline m_n^0(z_2)\)}\right]+o_p(1),
\end{align}
and
\begin{align}\label{inp3}
&\frac{b_{n}(z_1)b_{n}(z_2)}{pn}\sum_{k=1}^n\rtr\left(\bGa_n^*\re_k\bA_k^{-1}(z_1)\bGa_n\right)\circ\left(\bGa_n^*\re_k\bA_k^{-1}(z_2)\bGa_n\right)\notag\\
=&\frac{b_{n}(z_1)b_{n}(z_2)}p\rtr\left(\bGa_n^*\bbH_n^{-1}(z_1)\bGa_n\right)\circ\left(\bGa_n^*\bbH_n^{-1}(z_2)\bGa_n\right)+o_p(1)\notag\\
=&\frac{1}{c_{n2}^2}-\frac1{c_{n2}^2p}\rtr\left[{\rm diag}\(\left(\underline m_n^0(z_1)\bGa_n^*\bGa_n+\bI_p\right)^{-1}+\left(\underline m_n^0(z_2)\bGa_n^*\bGa_n+\bI_p\right)^{-1}\)\right]\\
&+\frac{1}{c_{n2}^2p}\rtr\(\left(\underline m_n^0(z_1)\bGa_n^*\bGa_n+\bI_p\right)^{-1}\circ\left(\underline m_n^0(z_2)\bGa_n^*\bGa_n+\bI_p\right)^{-1}\)+o_p(1)\notag
\\=&\frac{\frac{1}{p}\rtr\left(\bGa_n^*\(\underline m(z_1)\bSig_n+\bI_p\)^{-1}\bGa_n\right)\circ\left(\bGa_n^*\(\underline m(z_2)\bSig_n+\bI_p\)^{-1}\bGa_n\right)}{\(z_1-c_{1n}\int \frac t{1+t\underline m_n^0(z_1)}dH_n(t)\)\(z_2-c_{1n}\int \frac t{1+t\underline m_n^0(z_2)}dH_n(t)\)}+o_p(1),
\end{align}
where we used the fact that from \eqref{bnlimit}
$
	\frac{n}{\nu_{np}}b_n(z)=-z\underline m_n^0(z)+o(1).
$

On the other hand, denote $\breve\bA_k(z)=\frac1{\nu_{np}}\sum_{j<k}\bx_j\bx_j^*+\frac1{\nu_{np}}\sum_{j>k}\breve\bx_j\breve\bx_j^*-z\bI_p,$  where $\breve\bx_j$ is an i.i.d. copy of $\bx_j, j = k+1,\cdots,n.$ Thus, we deduce
\begin{align*}
&\frac{b_{n}(z_1)b_{n}(z_2)}p\re_k\left[z_1\rtr\left(\bA_k^{-1}(z_1)\bSig_n\right)-z_2\rtr\left(\breve\bA_k^{-1}(z_2)\bSig_n\right)\right]\\
=&\frac{b_{n}(z_1)b_{n}(z_2)}p\re_k\rtr\left[\bA_k^{-1}(z_1)\left(z_1\breve\bA_k(z_2)-z_2\bA_k(z_1)\right)\breve\bA_k^{-1}(z_2)\bSig_n\right]\\
=&\frac{b_{n}(z_1)b_{n}(z_2)}{p\nu_{np}}\re_k\rtr\bigg[\bA_k^{-1}(z_1)\bigg((z_1-z_2)\sum_{j=1}^{k-1}\bx_j\bx_j^*+z_1\sum_{j=k+1}^n\breve\bx_j\breve\bx_j^*\\
&\qquad\qquad\qquad\qquad\qquad\qquad-z_2\sum_{j=k+1}^n\bx_j\bx_j^*\bigg)\breve\bA_k^{-1}(z_2)\bSig_n\bigg]\\
=&\frac{kb_{n}^2(z_1)b_{n}^2(z_2)(z_1-z_2)}{p\nu_{np}}\re_k\rtr\left(\bA_{k}^{-1}(z_1)\bSig_n\breve\bA_{k}^{-1}(z_2)\bSig_n\right)+\left(z_1b_n(z_2)
-z_2b_n(z_1)\right)\\
&\qquad\cdot\frac{b_{n}(z_1)b_{n}(z_2)(n-k)}{p\nu_{np}}\re_k\rtr\(\bA_k^{-1}(z_1)\bSig_n\breve\bA_{k}^{-1}(z_2)\bSig_n\)+o_p(1)\\
=&z_1z_2\left(\underline m_n^0(z_1)-\underline m_n^0(z_2)
\right)\left[1-\frac{k}{n}\left(1-\frac{\underline m_n^0(z_1)\underline m_n^0(z_2)\(z_1-z_2\)}{c_{n2}\(\underline m_n^0(z_1)-\underline m_n^0(z_2)\)}\right)\right]\textcolor{red}{}\\
&\qquad\cdot\frac{b_{n}(z_1)b_{n}(z_2)}{p}\re_k\rtr\left(\bA_{k}^{-1}(z_1)\bSig_n\breve\bA_{k}^{-1}(z_2)\bSig_n\right)+o_p(1).
\end{align*}
By combining \eqref{con1} with the aforementioned equality, it becomes evident
\begin{align*}
&\frac{b_{n}(z_1)b_{n}(z_2)}{pn}\sum_{k=1}^n\rtr\(\re_k\bA_k^{-1}(z_1)\bSig_n\re_k\bA_k^{-1}(z_2)\bSig_n\)\\
=&\frac1{c_{n1}c_{n2}}\int_0^1\frac{1-\frac{\underline m_n^0(z_1)\underline m_n^0(z_2)\(z_1-z_2\)}{c_{n2}\(\underline m_n^0(z_1)-\underline m_n^0(z_2)\)}}{1-t\left(1-\frac{\underline m_n^0(z_1)\underline m_n^0(z_2)\(z_1-z_2\)}{c_{n2}\(\underline m_n^0(z_1)-\underline m_n^0(z_2)\)}\right)}dt+o_p(1)\textcolor{red}{}\\
=&\frac1{c_{n1}c_{n2}}\log\(\frac{c_{n2}\(\underline m_n^0(z_1)-\underline m_n^0(z_2)\)}{\underline m_n^0(z_1)\underline m_n^0(z_2)\(z_1-z_2\)}\)+o_p(1).\textcolor{red}{}
\end{align*}
Combining \eqref{inp}, \eqref{inp1}, \eqref{inp2}, \eqref{inp3} with the preceding equation, we observe
\begin{align*}
	\Phi(z_1,z_2)\xrightarrow{p}&\lim_{p,n\to\infty}\frac{1+\alpha_y}{c_{n1}c_{n2}}\left[\frac{\(\underline m_n^0(z_1)\)'\(\underline m_n^0(z_2)\)'}{\(\underline m_n^0(z_1)-\underline m_n^0(z_2)\)^2}-\frac1{\(z_1-z_2\)^2}\right]\textcolor{red}{}\\
	&+\Delta_yt(z_1,z_2)\lim_{p,n\to\infty}\frac{\(\underline m_n^0(z_1)\)'\(\underline m_n^0(z_2)\)'}{c_{n2}^2}.
\end{align*}
Our subsequent aim is to determine the limit of $$\frac{1}{c_{n1}c_{n2}}\left[\frac{\(\underline m_n^0(z_1)\)'\(\underline m_n^0(z_2)\)'\(z_1-z_2\)^2}{\(\underline m_n^0(z_1)-\underline m_n^0(z_2)\)^2}-1\right],\ {\rm and}\ \frac{\(\underline m_n^0(z_1)\)'\(\underline m_n^0(z_2)\)'}{c_{n2}^2}.$$
Using
\begin{align}\label{bnlim}
	b_n(z)\to\frac z{z-c_{1}\int\frac t{1+t\underline m(z)}dH(t)}=-\frac{z\underline m(z)}{c_2},
\end{align}
we find
\begin{align*}
	\underline m'(z)=&\(\frac{c_2}{c_1\int\frac t{1+t\underline m(z)}dH(t)-z}\)'=\frac{c_2}{g^2(z)-c_1c_2h(z,z)},
\end{align*}
and
\begin{align*}
&\underline m(z_1)-\underline m(z_2)=\frac{c_2\(z_1-z_2\)}{g(z_1)g(z_2)-c_1c_2h(z_1,z_2)}.
\end{align*}
From the above, we can obtain
\begin{align*}
	\frac{\underline m'(z_1)\underline m'(z_2)\(z_1-z_2\)^2}{\(\underline m(z_1)-\underline m(z_2)\)^2}-1
	=&\frac{\left[{g(z_1)g(z_2)-c_1c_2h(z_1,z_2)}\right]^2}{\left[g^2(z_1)-c_1c_2h(z_1,z_1)\right]\left[g^2(z_2)-c_1c_2h(z_2,z_2)\right]}-1\\
	=&\frac{c_1^2c_2^2\(h^2(z_1,z_2)-h(z_1,z_1)h(z_2,z_2)\)+c_1c_2g^2(z_1)h(z_2,z_2)}{\left[g^2(z_1)-c_1c_2h(z_1,z_1)\right]\left[g^2(z_2)-c_1c_2h(z_2,z_2)\right]}\\
	&+\frac{c_1c_2h(z_1,z_1)g^2(z_2)-2c_1c_2g(z_1)g(z_2)h(z_1,z_2)}{\left[g^2(z_1)-c_1c_2h(z_1,z_1)\right]\left[g^2(z_2)-c_1c_2h(z_2,z_2)\right]}.
\end{align*}
This implies
\begin{align*}
	&\frac{1}{c_{n1}c_{n2}}\left[\frac{\(\underline m_n^0(z_1)\)'\(\underline m_n^0(z_2)\)'\(z_1-z_2\)^2}{\(\underline m_n^0(z_1)-\underline m_n^0(z_2)\)^2}-1\right]\\
	\to&\frac{c_1c_2\(h^2(z_1,z_2)-h(z_1,z_1)h(z_2,z_2)\)+g^2(z_1)h(z_2,z_2)}{\left[g^2(z_1)-c_1c_2h(z_1,z_1)\right]\left[g^2(z_2)-c_1c_2h(z_2,z_2)\right]}\\
	&+\frac{h(z_1,z_1)g^2(z_2)-2g(z_1)g(z_2)h(z_1,z_2)}{\left[g^2(z_1)-c_1c_2h(z_1,z_1)\right]\left[g^2(z_2)-c_1c_2h(z_2,z_2)\right]}.
\end{align*}
Moreover, it is obvious
\begin{align*}
	\frac{\(\underline m_n^0(z_1)\)'\(\underline m_n^0(z_2)\)'}{c_{n2}^2}\to\frac{1}{\(g^2(z_1)-c_1c_2h(z_1,z_1)\)\(g^2(z_2)-c_1c_2h(z_2,z_2)\)}.
\end{align*}
Thus, we conclude
\begin{align*}
	&\Phi(z_1,z_2)\xrightarrow{p}\frac{1+\alpha_y}{\(z_1-z_2\)^2}\Bigg[\frac{c_1c_2\(h^2(z_1,z_2)-h(z_1,z_1)h(z_2,z_2)\)}{\left[g^2(z_1)-c_1c_2h(z_1,z_1)\right]\left[g^2(z_2)-c_1c_2h(z_2,z_2)\right]}\\
	&\qquad\qquad\quad+\frac{g^2(z_1)h(z_2,z_2)+h(z_1,z_1)g^2(z_2)-2g(z_1)g(z_2)h(z_1,z_2)}{\left[g^2(z_1)-c_1c_2h(z_1,z_1)\right]\left[g^2(z_2)-c_1c_2h(z_2,z_2)\right]}\Bigg]\\
	&\qquad\qquad\quad+\frac{\Delta_yt(z_1,z_2)}{\(g^2(z_1)-c_1c_2h(z_1,z_1)\)\(g^2(z_2)-c_1c_2h(z_2,z_2)\)}\\
	=&\begin{cases}
		\int\frac{\(1+\alpha_y\)t^2}{\(z_1-t\)^2\(z_2-t\)^2}dH(t)+\Delta_y\lim\frac{\rtr\(\bGa_n^*\(z_1\bI_p-\bSig_n\)^{-2}\bGa_n\)\circ\(\bGa_n^*\(z_2\bI_p-\bSig_n\)^{-2}\bGa_n\)}p,\ c_1=0,\\
		\frac{\(1+\alpha_y\)\int t^2dH(t)+\Delta_y\lim\frac1p\rtr\(\bGa_n^*\bGa_n\)\circ\(\bGa_n^*\bGa_n\)}{\(c_1\int tdH(t)-z_1\)^2\(c_1\int tdH(t)-z_2\)^2},\quad c_2=0,\\
		\frac{\(1+\alpha_y\)}{c_{1}c_{2}}\left[\frac{\underline m'(z_1)\underline m'(z_2)}{\(\underline m(z_1)-\underline m(z_2)\)^2}-\frac1{\(z_1-z_2\)^2}\right]+\frac{\Delta_y}{c_{2}^2}t(z_1,z_2)\underline m'(z_1)\underline m'(z_2),\quad c_1,c_2\neq 0.
	\end{cases}
\end{align*}

\paragraph{Convergence of $M_{n2}(z)$}
In this part, our objective is to determine the asymptotic behavior of  $M_{n2} (z)$ for $z\in\mathcal C_n$. To start, we observe that, similar to the comparable scenario, it can be easily verified that
$
	\sup_{z\in\mathcal C_n}|\re m_n(z)-m(z)|\to0,\ \sup_{z\in\mathcal C_n}\|\bH_n^{-1}(z)\|<C_0,
$
where $C_0$ is a constant. Moreover, we have $|b_n(z)|\le C_1$ and for $q\ge1,$
$$
	\max\{\re\|\bA^{-1}(z)\|^q,\re\|\bA_k^{-1}(z)\|^q,\re\|\bA_{kj}^{-1}(z)\|^q,\re|\beta_k(z)|^q\}\le C_q.
$$   It is important to note that, both here and in the proof of tightness, we essentially necessitate a high probability bound: $|\bS_n^0| \leq (1 + \sqrt{p/n})^2 |\bSig_n|$ under all three scenarios with high probability. This can be attained by following the entire procedure outlined in \cite{yin1988limit}, with a distinct truncation of variables at $\delta_n\sqrt[4]{pn}$. We omit the repetitive details of this process.

Returning to our primary objective, let's define
$a_n(z,\bM)=\frac{\nu_{np}}{\sqrt{pn}}[\re\rtr\(\bA^{-1}(z)\bM\)-\rtr\(\bH_n^{-1}(z)\bM\)\,$ where $\bM$ is a nonrandom matrix.
Note that
\begin{align*}
M_{n2}(z)=	&\frac{\nu_{np}}{\sqrt{pn}}\left[p\re\(m_n(z)\)-\rtr\(\bH_n^{-1}(z)\)\right]+\frac{\nu_{np}}{\sqrt{pn}}\left[\rtr\(\bH_n^{-1}(z)\)-pm_n^0(z)\right],
\end{align*}
and
\begin{align*}
&\frac{\nu_{np}}{\sqrt{pn}}\left[\rtr\(\bH_n^{-1}(z)\)-pm_n^0(z)\right]\\
=&-\frac{p\nu_{np}}{\sqrt{pn}}\(z\underline m_n^0(z)+\frac n{\nu_{np}}b_n(z)\)\int \frac{z^{-2}tdH_n(t)}{\(1+t\underline m_n^0(z)\)^2}+o(1)\\
=&\frac{\nu_{np}}{\sqrt{pn}}\Big(\rtr\(\bH_n^{-1}(z)\)- pm_n^0(z)\Big)\int \frac{c_{n1}z^{-1}tdH_n(t)}{\(1+t\underline m_n^0(z)\)^2}\\
&+\frac{\sqrt{pn}}{\nu_{np}}b_n(z)\Big(\re\rtr\(\bA^{-1}(z)\bSig_n\)-\rtr\(\bH_n^{-1}(z)\bSig_n\)\Big)\int \frac{z^{-2}tdH_n(t)}{\(1+t\underline m_n^0(z)\)^2}+o(1)\\
=&-\frac{\frac{\nu_{np}}{\sqrt{pn}}\Big(\re\(\rtr\bA^{-1}(z)\bSig_n\)-\rtr\(\bH_n^{-1}(z)\bSig_n\)\Big)\int \frac{c_{n1}z^{-1}\underline m_n^0(z)tdH_n(t)}{\(1+t\underline m_n^0(z)\)^2}}{1-c_{n1}z^{-1}\int \frac{tdH_n(t)}{\(1+t\underline m_n^0(z)\)^2}}+o(1).
\end{align*}
These imply
\begin{align}\label{m2limit}
M_{n2}(z)=	&a_n(z,\bI_p)-\frac{\int \frac{c_{n1}z^{-1}\underline m_n^0(z)tdH_n(t)}{\(1+t\underline m_n^0(z)\)^2}}{1-c_{n1}z^{-1}\int \frac{tdH_n(t)}{\(1+t\underline m_n^0(z)\)^2}}a_n(z,\bSig_n)+o(1).
\end{align}
 Hence, it suffices to derive the limit of $a_n(z,\bM)$ for $z\in\mathcal C_n$. 

To this end, we rewrite
\begin{align*}
&a_n(z,\bM)\\
=&-\frac{1}{\sqrt{pn}}\sum_{k=1}^n\re\beta_k(z)\(\bx_k^*\bA_k^{-1}(z)\bM\bH_n^{-1}(z)\bx_k-\re\rtr\(\bA_k^{-1}(z)\bM\bH_n^{-1}(z)\bSig_n\)\)\\
&-\frac{1}{\sqrt{pn}\nu_{np}}\sum_{k=1}^n\re\beta_k(z)\re\(\beta_k(z)\bx_k^*\bA_k^{-1}(z)\bM\bH_n^{-1}(z)\bSig_n\bA_k^{-1}(z)\bx_k\)\\
&-\frac{1}{\sqrt{pn}}\sum_{k=1}^n\(\re\beta_k(z)-b_n(z)\)\re\rtr\(\bA^{-1}(z)\bM
\bH_n^{-1}(z)\bSig_n\)\\
			\triangleq&\mathcal J_1(z)+\mathcal J_2(z)+\mathcal J_3(z).
\end{align*}
For a nonrandom $p \times p$ matrix $\bB$, employing the martingale difference decomposition method, we obtain
\begin{align*}
&\re\left|\rtr\(\bA_k^{-1}(z)\bB\bSig_n\)-\re\rtr\(\bA_k^{-1}(z)\bB\bSig_n\)\right|^q\le C_q\delta_n^{2q-4}n^{q/2-1},
\end{align*}
and
\begin{align*}
\re\left|\rtr\(\bA_k^{-1}(z)\bB\bSig_n\)-\re\rtr\(\bA^{-1}(z)\bB\bSig_n\)\right|^q
	\le&C_q\delta_n^{2q-4}n^{q/2-1}.
\end{align*}
Moreover, it follows from Lemma \ref{lep1} that
\begin{align*}
\re\left|\bx_k^*\bA_k^{-1}(z)\bB\bx_k-\re\rtr\(\bA_k^{-1}(z)\bB\bSig_n\)\right|^q
	\le&C_q\delta_n^{2q-4}p^{q/2}n^{q/2-1},
\end{align*}
and
\begin{align}\label{xi1}
	\re\left|\bx_k^*\bA_k^{-1}(z)\bB\bx_k-\re\rtr\(\bA^{-1}(z)\bB\bSig_n\)\right|^q
	\le&C_q\delta_n^{2q-4}p^{q/2}n^{q/2-1}.
\end{align}
Denoting $\xi_k(z)=\bx_k^*\bA_k^{-1}(z)\bx_k-\re\rtr\(\bA^{-1}(z)\bSig_n\)$, it yields
\begin{align}\label{beta}
	\beta_k(z)=b_n(z)-\frac1{\nu_{np}}\beta_k(z)b_n(z)\xi_k(z).
\end{align}
We will now examine $\mathcal J_k(z)$ for $k=1,2,3$. Focusing on $\mathcal J_1(z)$ and employing the aforementioned inequalities, we observe:
\begin{align*}
	&\mathcal J_1(z)\\
	=&\frac{b_n^2(z)}{\sqrt{pn}\nu_{np}}\sum_{k=1}^n\re\xi_k(z)\(\bx_k^*\bA_k^{-1}(z)\bM\bH_n^{-1}(z)\bx_k-\re\rtr\(\bA_k^{-1}(z)\bM\bH_n^{-1}(z)\bSig_n\)\)\\
	&-\frac{b_n^2}{\sqrt{pn}\nu_{np}^{2}}\sum_{k=1}^n\re\beta_k(z)\xi_k^2(z)\(\bx_k^*\bA_k^{-1}(z)\bM\bH_n^{-1}(z)\bx_k-\re\rtr\(\bA_k^{-1}(z)\bM\bH_n^{-1}(z)\bSig_n\)\)\\
	=&\frac{b_n^2(z)}{\sqrt{pn}\nu_{np}}\sum_{k=1}^n\re\xi_k(z)\(\bx_k^*\bA_k^{-1}(z)\bM\bH_n^{-1}(z)\bx_k-\re\rtr\(\bA_k^{-1}(z)\bM\bH_n^{-1}(z)\bSig_n\)\)+o(1)\\
	=&\frac{b_n^2(z)}{\sqrt{pn}\nu_{np}}\sum_{k=1}^n\re\varepsilon_k(z)\(\bx_k^*\bA_k^{-1}(z)\bM\bH_n^{-1}(z)\bx_k-\rtr\(\bA_k^{-1}(z)\bM\bH_n^{-1}(z)\bSig_n\)\)+o(1)\\
	=&\(1+\alpha_y\)\frac{b_n^2(z)}{\sqrt{pn}\nu_{np}}\sum_{k=1}^n\re\rtr\(\bA_k^{-1}(z)\bM\bH_n^{-1}(z)\bSig_n\bA_k^{-1}(z)\bSig_n\)\\
	&+\Delta_y\frac{b_n^2(z)}{\sqrt{pn}\nu_{np}}\sum_{k=1}^n\re\rtr\(\bGa_n^*\bA_k^{-1}(z)\bGa_n\)\circ\(\bGa_n^*\bA_k^{-1}(z)\bM\bH_n^{-1}(z)\bGa_n^*\)+o(1)\\
	=&\(1+\alpha_y\)\frac{ \sqrt {c_{n1}c_{n2}}b_n^2(z)}{p}\re\rtr\(\bA^{-1}(z)\bM\bH_n^{-1}(z)\bSig_n\bA^{-1}(z)\bSig_n\)\\
	&+\Delta_y\frac{ \sqrt {c_{n1}c_{n2}}b_n^2(z)}{p}\re\rtr\(\bGa_n^*\bA^{-1}(z)\bGa_n\)\circ\(\bGa_n^*\bA^{-1}(z)\bM\bH_n^{-1}(z)\bGa_n^*\)+o(1).
\end{align*}
Likewise, we can conclude that
\begin{align*}
	\mathcal J_2(z)=&-\frac{ b_n^2(z)}{\sqrt{pn}\nu_{np}}\sum_{k=1}^n\re\rtr\(\bA_k^{-1}(z)\bM\bH_n^{-1}(z)\bSig_n\bA_k^{-1}(z)\bSig_n\)+o(1)\\
	=&-\frac{ \sqrt {c_{n1}c_{n2}} b_n^2(z)}{p}\re\rtr\(\bA^{-1}(z)\bM\bH_n^{-1}(z)\bSig_n\bA^{-1}(z)\bSig_n\)+o(1),
\end{align*}
and
\begin{align*}
\mathcal J_3(z)=&\frac{b_n^2(z)}{\sqrt{pn}\nu_{np}}\sum_{k=1}^n\re\(\xi_k(z)\)\re\rtr\(\bA^{-1}(z)\bM\bH_n^{-1}(z)\bSig_n\)\\
&-\frac{b_n^3(z)}{\sqrt{pn}\nu_{np}^{2}}\sum_{k=1}^n\re\(\xi_k^2(z)\)\re\rtr\(\bA^{-1}(z)\bM\bH_n^{-1}(z)\bSig_n\)+o(1)\\
=&\frac{b_n^2(z)}{\sqrt{pn}\nu_{np}^2}\sum_{k=1}^n\re\(\beta_k(z)\bx_k^*\bA_k^{-1}(z)\bSig_n\bA_k^{-1}(z)\bx_k\)\re\rtr\(\bA^{-1}(z)\bM\bH_n^{-1}(z)\bSig_n\)\\
&-\frac{b_n^3(z)}{\sqrt{pn}\nu_{np}^{2}}\sum_{k=1}^n\re\(\varepsilon_k^2(z)\)\re\rtr\(\bA^{-1}(z)\bM\bH_n^{-1}(z)\bSig_n\)+o(1)\\
=&-\alpha_y\frac{ c_{n1}^{3/2}c_{n2}^{1/2}b_n^3(z)}{p^{2}}\re\rtr\(\bA^{-1}(z)\bSig_n\)^2\re\rtr\(\bA^{-1}(z)\bM\bH_n^{-1}(z)\bSig_n\)\\
&-\Delta_y\frac{ c_{n1}^{3/2}c_{n2}^{1/2}b_n^3(z)}{p^{2}}\re\rtr\(\bGa_n^*\bA^{-1}(z)\bGa_n\)\circ\(\bGa_n^*\bA^{-1}(z)\bGa_n^*\)\\
&\qquad\qquad\qquad\qquad\cdot\re\rtr\(\bA^{-1}(z)\bM\bH_n^{-1}(z)\bSig_n\)+o(1).
\end{align*}
Hence, $a_n(z,\bM)$ can be represented as
\begin{align*}
	a_n(z,&\bM)=\alpha_y\frac{ \sqrt {c_{n1}c_{n2}} b_n^2(z)}{p}\re\rtr\(\bA^{-1}(z)\bM\bH_n^{-1}(z)\bSig_n\bA^{-1}(z)\bSig_n\)\\
	&+\Delta_y\frac{ \sqrt {c_{n1}c_{n2}} b_n^2(z)}{p}\re\rtr\(\bGa_n^*\bA^{-1}(z)\bGa_n\)\circ\(\bGa_n^*\bA^{-1}(z)\bM\bH_n^{-1}(z)\bGa_n\)\\
	&-\alpha_y\frac{ c_{n1}^{3/2}c_{n2}^{1/2}b_n^3(z)}{p^{2}}\re\rtr\(\bA^{-1}(z)\bSig_n\)^2\re\rtr\(\bA^{-1}(z)\bM\bH_n^{-1}(z)\bSig_n\)\\
	&-\Delta_y\frac{ c_{n1}^{3/2}c_{n2}^{1/2}b_n^3(z)}{p^{2}}\re\rtr\(\bGa_n^*\bA^{-1}(z)\bGa_n\)\circ\(\bGa_n^*\bA^{-1}(z)\bGa_n\)\\
	&\qquad\qquad\qquad\qquad\cdot\re\rtr\(\bA^{-1}(z)\bM\bH_n^{-1}(z)\bSig_n\)+o(1).
\end{align*}
Recall that
\begin{align*}
	\bA^{-1}(z)-\bH_n^{-1}(z)
=&-\frac{b_n(z)}{\nu_{np}}\sum_{k=1}^n\bH_n^{-1}(z)\(\bx_k\bx_k^*-\bSig_n\)\bA_k^{-1}(z)\\
	&-\frac1{\nu_{np}}\sum_{k=1}^n\(\beta_k(z)-b_n(z)\)\bH_n^{-1}(z)\bx_k\bx_k^*\bA_k^{-1}(z)\\
	&+\frac {b_n(z)}{\nu_{np}}\sum_{k=1}^n\bH_n^{-1}(z)\bSig_n\(\bA^{-1}(z)-\bA_k^{-1}(z)\)\\
	\triangleq&\bC_1(z)+\bC_2(z)+\bC_3(z).
\end{align*}
Consider ${\bbM}$ as a $p\times p$ matrix with bounded moments of the spectral norm. By applying Lemma \ref{lep1}, \eqref{beta}, and \eqref{xi1}, we derive
\begin{align}\label{al16}
	&\re|\beta_{k}(z)-b_n(z)|^2
	\le\frac C{\nu_{np}^2}\re|\beta_k(z)\xi_k(z)|^2\le\frac {C\delta_n^2}{\sqrt{\nu_{np}}}
\end{align}
which implies that
\begin{align*}
	\re\bigg|\rtr\big({\bC_2}(z)&{\bbM}\big)\bigg|
	\le\frac {Cp}{\nu_{np}}\sum_{k=1}^n\re^{1/2}\left|\beta_{k}(z)-b_n(z)\right|^2
	\le\frac {C\delta_n^2pn}{\nu_{np}^{5/4}}.
\end{align*}
Using Lemma \ref{lep1}, we establish that
$
	\re\left|\rtr\left({\bC_3}(z){\bbM}\right)\right|
	\le C.
$
Subsequently, assuming ${\bbM}$ is a nonrandom matrix, it becomes evident that
$
	\re\left|\rtr\left({\bC_1}(z){\bbM}\right)\right|
	\le \frac{Cn\sqrt p}{\nu_{np}}.
$
Combining these inequalities with the previously mentioned three, we deduce that
\begin{align}\label{an1}
\begin{split}
\frac1{p}\re\rtr\big(\bA^{-1}(z)\bM\bH_n^{-1}(z)&\bSig_n\big)=\frac1{p}\re\rtr\(\bH_n^{-1}(z)\bM\bH_n^{-1}(z)\bSig_n\)+o(1),
\end{split}
\end{align}
and
\begin{align*}
	&\frac1{p}\re\rtr\(\bSig_n\bA^{-1}(z)\bSig_n\bA^{-1}(z){\bf M}\bbM\)\\
=&\frac1{p}\re\rtr\(\bSig_n\bA^{-1}(z)\bSig_n\bH_n^{-1}(z){\bf M}\bbM\)+\frac1{p}\re\rtr\(\bSig_n\bA^{-1}(z)\bSig_n\bC_1(z){\bf M}\bbM\)+o(1)\\
=&\frac1{p}\rtr\(\bSig_n\bH_n^{-1}(z)\bSig_n\bH_n^{-1}(z){\bf M}\bbM\)+\frac1{p}\re\rtr\(\bSig_n\bA^{-1}(z)\bSig_n\bC_1(z){\bf M}\bbM\)+o(1).
\end{align*}
We are now ready to determine the limit of  $p^{-1}\re\rtr\(\bSig_n\bA^{-1}(z)\bSig_n\bC_1(z){\bf M}\bbM\)$. Upon calculation, it implies that
\begin{align*}
&\frac1{p}\re\rtr\(\bSig_n\bA^{-1}(z)\bSig_n\bC_1(z){\bf M}\bbM\)\\
=&-\frac{b_n(z)}{p\nu_{np}}\sum_{k=1}^n\re\rtr\(\bSig_n\bA^{-1}(z)\bSig_n\bH_n^{-1}(z)\(\bx_k\bx_k^*-\bSig_n\)\bA_k^{-1}(z){\bf M}\bbM\)\\
=
&\frac{b_n^2(z)}{p\nu_{np}^2}\sum_{k=1}^n\re\(\bx_k^*\bA_k^{-1}(z)\bSig_n\bH_n^{-1}(z)\bx_k\bx_k^*\bA_k^{-1}(z){\bf M}\bbM\bSig_n\bA_k^{-1}(z)\bx_k\)+o(1)\\
=
&\frac{b_n^2(z)}{p\nu_{np}^2}\sum_{k=1}^n\re\rtr\(\bA_k^{-1}(z)\bSig_n\bH_n^{-1}(z)\bSig_n\)\rtr\(\bA_k^{-1}(z){\bf M}\bbM\bSig_n\bA_k^{-1}(z)\bSig_n\)+o(1)\\
=
&\frac{nb_n^2(z)}{p\nu_{np}^2}\re\rtr\(\bA^{-1}(z)\bSig_n\bH_n^{-1}(z)\bSig_n\)\re\rtr\(\bA^{-1}(z){\bf M}\bbM\bSig_n\bA^{-1}(z)\bSig_n\)+o(1)\\
=
&\frac{nb_n^2(z)}{p\nu_{np}^2}\re\rtr\(\bH_n^{-1}(z)\bSig_n\bH_n^{-1}(z)\bSig_n\)\re\rtr\(\bA^{-1}(z){\bf M}\bbM\bSig_n\bA^{-1}(z)\bSig_n\)+o(1)\\
=
&-\frac{b_n(z)}{\nu_{np}}\int\frac{z^{-1}\underline m_n^0(z)t^2}{\(1+t\underline m_n^0(z)\)^2}\re\rtr\(\bSig_n\bA^{-1}(z)\bSig_n\bA^{-1}(z){\bf M}\bbM\)+o(1).
\end{align*}
Combining with the above equalities, we observe that
\begin{align*}
\frac1{p}\re\rtr\big(\bSig_n\bA^{-1}(z)	&\bSig_n\bA^{-1}(z){\bf M}\bbM\big)
	=\frac{p^{-1}\rtr\(\bSig_n\bH_n^{-1}(z)\bSig_n\bH_n^{-1}(z){\bf M}\bbM\)}{1+c_{n1}{b_n(z)}\int\frac{z^{-1}\underline m_n^0(z)t^2}{\(1+t\underline m_n^0(z)\)^2}dH_n(t)}+o(1)\\
	=&\frac{p^{-1}\rtr\(\bSig_n\bH_n^{-1}(z)\bSig_n\bH_n^{-1}(z){\bf M}\bbM\)}{1-\frac{c_{n1}}{c_{n2}}\int\frac{\(\underline m_n^0(z)t\)^2}{\(1+t\underline m_n^0(z)\)^2}dH_n(t)}+o(1)
	.
\end{align*}
Extending the same argument, it follows that
\begin{align*}
	&\frac1p\re\rtr\(\bGa_n^*\bA^{-1}(z)\bGa_n\)\circ\(\bGa_n^*\bA^{-1}(z){\bf M}\bbM\bGa_n\)\\
	=&\frac1p\rtr\(\bGa_n^*\bH_n^{-1}(z)\bGa_n\)\circ\(\bGa_n^*\bH_n^{-1}(z){\bf M}\bbM\bGa_n\)+o(1).
\end{align*}
Combining \eqref{an1} with the three aforementioned equalities, we obtain
\begin{align*}
	a_n(z,&\bM)=\alpha_y\frac{ \sqrt {c_{n1}c_{n2}} \(z\underline m_n^0(z)\)^2}{pc_{n2}^2}\frac{\rtr\(\bSig_n\bH_n^{-1}(z)\bSig_n\bH_n^{-1}(z){\bf M}\bH_n^{-1}(z)\)}{\(1-\frac{c_{n1}}{c_{n2}}\int\frac{\(\underline m_n^0(z)t\)^2}{\(1+t\underline m_n^0(z)\)^2}dH_n(t)\)}\\
		&+\alpha_y\frac{ c_{n1}^{3/2}c_{n2}^{1/2}z\(\underline m_n^0(z)\)^3\int\frac{t^2}{\(1+t\underline m_n^0(z)\)^2}dH_n(t)}{c_{n2}^3p\(1-\frac{c_{n1}}{c_{n2}}\int\frac{\(\underline m_n^0(z)t\)^2}{\(1+t\underline m_n^0(z)\)^2}dH_n(t)\)}\rtr\(\bH_n^{-1}(z)\bM\bH_n^{-1}(z)\bSig_n\)\\
	&+\Delta_y\frac{ \sqrt {c_{n1}c_{n2}} \(z\underline m_n^0(z)\)^2}{pc_{n2}^2}\re\rtr\(\bGa_n^*\bH_n^{-1}(z)\bGa_n\)\circ\(\bGa_n^*\bH_n^{-1}(z)\bM\bH_n^{-1}(z)\bGa_n\)\\
	&+\Delta_y\frac{ c_{n1}^{3/2}c_{n2}^{1/2}\(z\underline m_n^0(z)\)^3}{p^{2}c_{n2}^3}\re\rtr\(\bGa_n^*\bH_n^{-1}(z)\bGa_n\)\circ\(\bGa_n^*\bH_n^{-1}(z)\bGa_n\)\\
	&\qquad\qquad\qquad\qquad\cdot\re\rtr\(\bH_n^{-1}(z)\bM\bH_n^{-1}(z)\bSig_n\)+o(1).
\end{align*}
By combining the previously mentioned equality, \eqref{bnlim}, and \eqref{m2limit}, we deduce
\begin{align*}
	M_{n2}(z)-
\(\frac{\alpha_y\sqrt {c_{n1}c_{n2}}\(\underline m_n^0(z)\)^3\int\frac{t^2}{\(1+t\underline m(z)\)^3}dH(t)}{c_{n2}^3\(1-\frac{c_{n1}}{c_{n2}}\int\frac{\(\underline m_n^0(z)t\)^2}{\(1+t\underline m_n^0(z)\)^2}dH(t)\)^2}+\frac{\Delta_y\sqrt {c_{n1}c_{n2}}\(\underline m_n^0(z)\)^3s(z) }{c_{n2}^3\(1-\frac{c_{n1}}{c_{n2}}\int\frac{\(\underline m_n^0(z)t\)^2}{\(1+t\underline m_n^0(z)\)^2}dH(t)\)}\)\to 0
\end{align*}

\paragraph{Tightness of $M_{n1}(z)$}
The demonstration of the tightness of $M_{n1}(z)$ indeed relies on the application of Theorem 12.3 from \cite{Billingsley68}. It is therefore sufficient to demonstrate
\begin{align*}
	\re\frac{\left|M_{n1}(z_1)-M_{n1}(z_2)\right|^2}{|z_1-z_2|^2}\le C, \quad z_1,z_2\in\mathcal C_n.
\end{align*}
Note that
\begin{align*}
&\frac{M_{n1}(z_1)-M_{n1}(z_2)}{z_1-z_2}\\
=&\frac{\nu_{np}}{\sqrt{pn}}\sum_{k=1}^n\(\re_k-\re_{k-1}\)\bigg[\rtr\(\bA^{-1}(z_1)\bA^{-1}(z_2)\)-\rtr\(\bA_k^{-1}(z_1)\bA_k^{-1}(z_2)\)\bigg]\\
=&\frac{1}{\sqrt{pn}\nu_{np}}\sum_{k=1}^n\(\re_k-\re_{k-1}\)\beta_k(z_1)\beta_k(z_2)\(\bx_k^*\bA_k^{-1}(z_1)\bA_k^{-1}(z_2)\bx_k\)^2\\
&-\frac{1}{\sqrt{pn}}\sum_{k=1}^n\(\re_k-\re_{k-1}\)\(\beta_k(z_1)\bx_k^*\bA_k^{-2}(z_1)\bA_k^{-1}(z_2)\bx_k\)\\
&-\frac{1}{\sqrt{pn}}\sum_{k=1}^n\(\re_k-\re_{k-1}\)\(\beta_k(z_2)\bx_k^*\bA_k^{-1}(z_1)\bA_k^{-2}(z_2)\bx_k\)\\
\triangleq&W_1(z_1,z_2)+W_2(z_1,z_2)+W_3(z_1,z_2).
\end{align*}
Obtaining
$
	\re\left|W_j(z_1,z_2)\right|^2\le C,\  j=1,2,3,
$
is straightforward following the same argument as in \cite{BaiS04C}, and as such, we omit this analogous derivation.
\paragraph{Completing the proof} Combining the aforementioned arguments, we conclude the proof of Lemma \ref{CLTST} and, consequently, the proof of Theorem \ref{cltthm}.

\subsection{Proof of Theorem \ref{identitytestf} and Theorem \ref{identitytestl}}
The proofs of these two theorems both hinge on the fundamental computation of a specific integral using the residue theorem. Consequently, we defer these proofs to Appendix \ref{cltappli}.

\section{Auxiliary lemmas}
Here, we outline and present the lemmas that form the backbone of our analysis.

\begin{lemma}[Corollary A.42 of \cite{BaiS10S}] \label{lemma:4}
Let ${\mathbf A}$ and ${\mathbf B}$ be two $p \times n$  matrices and denote the ESD of ${\mathbf S} = {\mathbf A}{{\mathbf A}^ * }$ and $ \widetilde{\mathbf S} = {\mathbf B}{{\mathbf B}^ * }$   by ${F^{\mathbf S}}$ and ${F^{\widetilde{\mathbf S}}}$, respectively.  Then,
$${L^4}\left({F^{\mathbf S}},{F^{\widetilde{\mathbf S}}}\right) \le \frac{2}{{{p^2}}}\left({\rm tr}\left({\mathbf A}{{\mathbf A}^ * } + {\mathbf B}{{\mathbf B}^ * }\right)\right)\left({\rm tr}\left[\left({\mathbf A} - {\mathbf B}\right){\left({\mathbf A} - {\mathbf B}\right)^ * }\right]\right),$$
 where $L\left(\cdot,\cdot \right)$  denotes the L\'{e}vy distance, that is,
 $$L\left({F^{\mathbf S}},{F^{\widetilde{\mathbf S}}}\right)=\inf\left\{\varepsilon:{F^{\mathbf S}\left(x-\varepsilon,y-\varepsilon\right)-\varepsilon}\le{F^{\widetilde{\mathbf S}}}\le{F^{\mathbf S}\left(x+\varepsilon,y+\varepsilon\right)+\varepsilon}\right\}.$$
\end{lemma}

\begin{lemma}[Theorem A.44 of \cite{BaiS10S}]\label{lemma:3}
Let ${\mathbf A}$ and ${\mathbf B}$ be $p \times n$ complex matrices. Then, $\left\| {{F^{{\mathbf A}{{\mathbf A}^ * }}} - {F^{{\mathbf B}{{\mathbf B}^ * }}}} \right\|_{KS} \le \frac{1}{p}{\rm rank}\left({\mathbf A} - {\mathbf B}\right).$
\end{lemma}

\begin{lemma}[Rosenthal's inequality ]\label{lemma:5}
Let $ \{{\mathbf  X}_k\} $ be a complex martingale difference sequence with respect to the increasing $\sigma$-field $\mathcal{F}_k$, and let $\re_k$ denote conditional expectation given $\mathcal{F}_k$. Then we have, for $q\ge2$,
$${\rm E}\left|\sum_k {\mathbf  X}_k\right|^{q} \leq C_q\left(\sum_k{\rm E}\left|{\mathbf  X}_k\right|^{q}+\re\left(\sum_k {\rm E}_{k-1}\left|{\mathbf  X}_k\right|^2\right)^{q/2}\right). $$
\end{lemma}

\begin{lemma}[Vitali's convergence theorem]\label{lemma:12}
Let $f_1,f_2,\cdots$ be analytic in $D$, a connected open set of $\mathbb C$, satisfying $\left|f_n\left(z\right)\right|\le M$ for every $n$ and $z$ in $D$, and $f_n\left(z\right)$ converges as $n \to \infty$ for each $z$ in a subset of $D$ having a limit point in $D$. Then there exists a function $f$ analytic in $D$ for which $f_n\left(z\right) \to f\left(z\right)$ and $f_n^{\prime} \to f^{\prime}\left(z\right)$ for all $z\in D$. Moreover, on any set bounded by a contour interior to $D$, the convergence is uniform and $\left\{f_n^{\prime}\left(z\right)\right\}$ is uniformly bounded.
\end{lemma}

\begin{lemma}[Lemma B.26 in \cite{BaiS10S}]\label{lep1}
	Let ${\bf A}$ be an $n\times n$ nonrandom matrix and $\bx=(x_1,\cdots,x_n)^T$ be a random vector of independent entries. Assume that $\re x_j=0$, $\re|x_j|^2=1$ and $\re|x_j|^{\ell}\le \nu_{\ell}$. Then for $q\ge2$,
	\begin{align*}
		&\re\left|\bx^*{\bf A}\bx-\rtr{\bf A}\right|^q
		\le C_q\left[\left(\nu_4\rtr{\bf A}{\bf A}^*\right)^{q/2}+\nu_{2q}\rtr\left({\bf A}{\bf A}^*\right)^{q/2}\right]
	\end{align*}
	where $C_q$ is a constant depending on $q$ only.
\end{lemma}
\begin{lemma}[CLT for martingale]\label{lep2}
	Suppose for each $n$ $Y_{n1},Y_{n2},\cdots,Y_{nr_n}$ is a real martingale difference sequence with respect to the increasing $\sigma$-field $\{\mathcal{F}_{nj}\}$ having second moments. If as $n\to\infty$,
	$(i): \sum_{j=1}^{r_n}\re(Y_{nj}^2|\mathcal{F}_{n,j-1})\xrightarrow{i.p.}\sigma^2,$
		 where $ \sigma^2$ is a positive constant, and for each $\varepsilon\ge0, $
		$(ii)\sum_{j=1}^{r_n}\re(Y_{nj}^2I(|Y_{nj}\ge\varepsilon|))\to0,$
		then $\sum_{j=1}^{r_n}Y_{nj}\xrightarrow{D}N(0,\sigma^2).
$
\end{lemma}

\section{Conclusion discussion}
In this paper, we introduce a universal approach to analyze the spectral properties of the sample covariance matrix, irrespective of the rates of divergence for both dimension  $p$ and sample size $n$. Our approach also addresses the issue of an unbounded spectral norm, which frequently arises in practical scenarios. We derive the LSD under three distinct scenarios and subsequently establish the CLT for the LSS. To demonstrate the efficacy of our theoretical results, we apply our approach to an identity test of population covariance matrix. We highlighting the dimension-proof property of the test statistics in this test. Our findings hold particular significance in practical data analysis, where both $p$ and $n$ are treated as constant values. By removing restrictions on their rates, our results become more reasonable and applicable. 

\appendix

\section{Details in the proof of Theorem \ref{th:1}}\label{Dtlsd}
Without loss of generality, in the proof of Theorem \ref{th:1}, we assume that $\re\(\bx_1\)=0$. In fact, let
$
\mathbf T_n=\frac{1}{\nu_{np}}\left(\mathbf{X}_n-{\rm E}{\bf X}_n\right)\left(\mathbf{X}_n-{\rm E}{\bf X}_n\right)^*.
$
By Lemma \ref{lemma:3}, it follows that, when $p\to\infty,$
$
\left\| {{F^{\mathbf S_n}} - {F^{\mathbf T_n}}} \right\|_{KS} \le \frac{1}{p}rank\left({\rm E}\mathbf X_n\right)= \frac1p\to 0,
$
where $\left\|f\right\|_{KS}=\sup_x\left|f(x)\right|$. 
Denote
\begin{align}\label{aa}
{ m}_n\left(z\right)=\frac{1}{p}{\rm tr}\left(\mathbf S_n-z\mathbf I_{p}\right)^{-1},
\end{align}
where $z=u+\upsilon i\in\mathbb{C}^+$. 

\subsection {Step 1}\label{se:1}
We begin by showing 
\begin{equation}\label{eq:6}{{ m}_n}\left(z\right) - {\rm E}{{m}_n}\left(z\right) \stackrel{a.s.}{\longrightarrow}0.\end{equation}
To establish this almost sure convergence, delve into the expression of the difference ${{m}_n}\left(z\right) - {\rm E}{{m}_n}\left(z\right)$ in terms of conditional expectations.  Let ${\re_k}\left( \cdot \right)$ denote the conditional expectation given $\left\{ \bx_{ 1},{\bx_{2}}, \cdots ,\bx_k \right\}$. Utilizing the Sherman-Morrison-Woodbury formula \begin{align}\label{inver}
	{\left({\mathbf A} + {\boldsymbol{\alpha}} {{\boldsymbol{\beta}}^*}\right)^{ - 1}} = {{\mathbf A}^{ - 1}} - \frac{{{{\mathbf A}^{ - 1}}{\boldsymbol{\alpha} } {{\boldsymbol{\beta}} ^ * }{{\mathbf A}^{ - 1}}}}{{1 + {{\boldsymbol{\beta}}^ * }{{\mathbf A}^{ - 1}}{\boldsymbol{\alpha} }}},
\end{align} for $\ell=0,1$, we have: 
\begin{align*}
&\frac1p\rtr\(\bA^{-1}(z)\bSig_n^{\ell}\)-\frac1p\re\rtr\(\bA^{-1}(z)\bSig_n^{\ell}\)\\
 =& \frac1p\sum_{k=1}^n\(\re_k-\re_{k-1}\)\left[\rtr\(\bA^{-1}(z)\bSig_n^{\ell}\)-\rtr\(\bA_k^{-1}(z)\bSig_n^{\ell}\)\right]\\
 =&-\frac{1}{{p\nu_{np}}}\sum_{k=1}^n\(\re_k-\re_{k-1}\)\beta_k(z)\bx_k^*\bA_k^{-1}(z)\bSig_n^{\ell}\bA_k^{-1}(z)\bx_k.
\end{align*}
Applying Lemma \ref{lemma:5} and considering condition $(c)$, we can bound the fourth moment:
\begin{align}\label{ran}
	\begin{split}
	 &{\rm E}{\left|\frac1p\rtr\(\bA^{-1}(z)\bSig_n^{\ell}\)-\frac1p\re\rtr\(\bA^{-1}(z)\bSig_n^{\ell}\)\right|^4}\\
	  \le& \frac{Cn^2}{p^4\nu_{np}^4}{\rm E}\left(\bx_1^*\bA_1^{-1}(z)\bSig_n^{\ell}\bA_1^{-1}(z)\bx_1\right)^4= O\left(p^{-2}{n^{ - 2}}\right).
	  \end{split}
\end{align}
Applying the Borel-Cantelli lemma and Chebyshev inequality to the fourth moment bound \eqref{ran}, we conclude  $${{ m}_n}\left(z\right) - {\rm E}{{ m}_n}\left(z\right)\stackrel{a.s.}{\longrightarrow}0.$$

\subsection{Step 2}\label{se:2}
Write
$
	\bH_n(z)=\frac n{\nu_{np}}b_n(z)\bSig_n-z\bI_p.
$
Then for any $t>0$,
\begin{align*}
	\left|\frac1{\frac n{\nu_{np}}b_n(z)t-z}\right|
	=&\left|\frac{\nu_{np}+\re\rtr\(\bA^{-1}(z)\bSig_n\)} {nt-z\(\nu_{np}+\re\rtr\(\bA^{-1}(z)\bSig_n\)\)}\right|
	\le \frac1v\(1+\frac1{\nu_{np}}\left|\re\rtr\(\bA^{-1}(z)\bSig_n\)\right|\)\\
	\le&\frac1v\(1+\frac {p\|\bSig_n\|}{v\nu_{np}}\)=\frac{C}{v^2}.
\end{align*}
This yields
\begin{align}\label{hbound}
	\left\|\bH_n^{-1}(z)\right\|\le \frac C{v^2}.
\end{align}

Note that from \eqref{inver}, 
\begin{align*}
&\bA^{-1}(z)-\bH_n^{-1}(z)=\bH_n^{-1}(z)\(\bH_n-\bA\)\bA^{-1}(z)\\
=&-\frac1{\nu_{np}}\sum_{k=1}^n\beta_k(z)\bH_n^{-1}(z)\bx_k\bx_k^*\bA_k^{-1}(z)+\frac n{\nu_{np}}b_n(z)\bH_n^{-1}(z)\bSig_n\bA^{-1}(z).
\end{align*}
Multiplying $\bSig_n^{\ell}, \ell=0,1$ on both sides of the above equality,
it follows that by simple calculation:
\begin{align}\label{meancom}
\begin{split}&\re\(m_n(z)\)-\frac1p\rtr\(\bH_n^{-1}(z)\)\\
=&-\frac1{p\nu_{np}}\sum_{k=1}^n\re\left[\(\beta_k(z)-b_n(z)\)\bx_k^*\bA_k^{-1}(z)\bH_n^{-1}(z)\bx_k\right]\\
&+\frac {b_n(z)}{p\nu_{np}}\sum_{k=1}^n\re\left[\rtr\(\bH_n^{-1}(z)\bSig_n\bA^{-1}(z)\)-\rtr\(\bH_n^{-1}(z)\bSig_n\bA_k^{-1}(z)\)\right]\\
\triangleq&\mathcal I_1(z)+\mathcal I_2(z).
\end{split}\end{align}
and
\begin{align*}
&\frac1n\re\rtr\(\bA^{-1}(z)\bSig_n\)-\frac1n\rtr\(\bH_n^{-1}(z)\bSig_n\)\notag\\
		=&-\frac1{n\nu_{np}}\sum_{k=1}^n\re\left[\(\beta_k(z)-b_n(z)\)\bx_k^*\bA_k^{-1}(z)\bSig_n\bH_n^{-1}(z)\bx_k\right]\\
		&+\frac {b_n(z)}{n\nu_{np}}\sum_{k=1}^n\re\left[\rtr\(\bH_n^{-1}(z)\bSig_n\bA^{-1}(z)\bSig_n\)-\rtr\(\bH_n^{-1}(z)\bSig_n\bA_k^{-1}(z)\bSig_n\)\right].\notag
	\end{align*}
By \eqref{inver} and condition $(c)$, it is evident
\begin{align}\label{al2}
	\left|\mathcal I_2(z)\right|\le& \frac{C}{p\nu_{np}^2}\sum_{k=1}^n\re^{1/2}\left|\bx_k^*\bA_k^{-1}(z)\bSig_n^{\ell}\bH_n^{-1}(z)\bSig_n\bA_k^{-1}(z)\bx_k\right|^2
	\le\frac{C}{\nu_{np}}\to0.
\end{align}
Utilizing the condition $(c)$, \eqref{ran}, and \eqref{inver}, we determine
\begin{align*}
	\re\left|\beta_k(z)-b_n(z)\right|^2
	\le&\frac C{\nu_{np}^2}\re\left|\bx_k^*\bA_k^{-1}(z)\bx_k-\re\rtr\(\bA^{-1}(z)\bSig_n\)\right|^2\\
\le&\frac C{\nu_{np}^2}\re\left|\rtr\(\bA_k^{-1}(z)\bSig_n\)-\rtr\(\bA^{-1}(z)\bSig_n\)\right|^2+\frac C{\nu_{np}}
\le \frac C{\nu_{np}}.
\end{align*}
Hence, it is apparent
\begin{align}\label{al1}
	\left|\mathcal I_1(z)\right|\le& \frac{C}{\nu_{np}}\sum_{k=1}^n\re^{1/2}\left|\beta_k(z)-b_n(z)\right|^2
	\le\frac{C}{\sqrt{\nu_{np}}}\to0.
\end{align}
Comparing (\ref{meancom}), \eqref{al2} with (\ref{al1}), we obtain
\begin{align}\label{nonran}
&\re\(m_n(z)\)-\frac1p\rtr\(\bH_n^{-1}(z)\) \to 0.
\end{align}
Applying a similar method to \eqref{nonran}, we establish the following result
\begin{align}
&\frac1n\re\rtr\(\bA^{-1}(z)\bSig_n\)-\frac1n\rtr\(\bH_n^{-1}(z)\bSig_n\) \to 0.\label{bn}
\end{align}
Recalling the definition of $\bH_n(z)$ and multiplying $\frac{n^2}{\nu_{np}^2}b_n(z)$ on both sides of \eqref{bn}, we obtain
\begin{align*}
	\frac{n}{\nu_{np}}-\frac{n}{\nu_{np}}b_n(z)-\frac{p}{\nu_{np}}-\frac{z}{\nu_{np}}\rtr\(\bH_n^{-1}(z)\)  \to 0.
\end{align*}
 Together with \eqref{nonran}, one has
\begin{align}\label{bnlimit}
	\frac{n}{\nu_{np}}b_n(z)=\frac{n}{\nu_{np}}-\frac{p}{\nu_{np}}\(1+z\re\(m_n(z)\)\)+o(1).
\end{align}
Substituting this into \eqref{nonran}, we see
\begin{align*}
	&\re\(m_n(z)\)-\int\frac1{\(\frac{n}{\nu_{np}}-\frac{p}{\nu_{np}}\(1+z\re\(m_n(z)\)\)\)t-z}dH_n(t)\to 0.
\end{align*}
For each fixed $z$, the sequence $\{\re(m_n(z))\}$ is bounded. Therefore, for any subsequence, they must converge to the same limit, denoted as $m(z)$. Then $m(z)$ should satisfy the equation \eqref{equa}, i.e.,
$
	m(z)=\int\frac1{\(c_2-c_1-c_1zm(z)\)t-z}dH(t).
$
Since $\Im(\re m_n(z)) > 0$, we conclude that $\Im(m(z)) \geq 0$. Obviously, it is impossible that $\Im(m(z)) = 0$ because the right-hand side of \eqref{equa} has a positive imaginary part.

\subsection{Step 3: Uniqueness of the solution of \eqref{equa}}
Let 
$\underline m(z)=-(c_2-c_1)/z+c_1m(z).$
Suppose we have two solutions $m_1(z), m_2(z) \in \mathbb{C}^{+}$ of the equation \eqref{equa}, namely
\begin{align*}
	&m_j(z)=\int\frac1{\(c_2-c_1-c_1zm_j(z)\)t-z}dH(t),\quad j=1,2.
\end{align*}
and
\begin{align}\label{uni1}
	&\underline m_j(z)=\frac{c_2}{c_1\int\frac t{1+\underline m_j(z)t}dH(t)-z},\quad j=1,2,
\end{align}
where $\underline m_j(z)=-(c_2-c_1)/z+c_1m_j(z)$. If $c_1c_2 = 0$, we assert that \eqref{equa} has at most one solution. When $c_2 = 0$, we obtain from
 \eqref{uni1}
\begin{align*}
	\underline m_1(z)=\underline m_2(z)=0\Rightarrow  m_1(z)= m_2(z)=-\frac1z.
\end{align*}
Thus the conclusion holds. Thus, the conclusion holds. Otherwise, the result can be inferred from
  $$m_j(z)=\int\frac1{tc_2-z}dH(t),\quad j=1,2.$$

We are now in a position to demonstrate that \eqref{equa} has at most one solution under $c_1c_2\neq 0$ and $\underline{m}_j(z)\in\mathbb{C}^+$. Note that
\begin{align*}
	&\underline m_1(z)-\underline m_2(z)
	=\frac{c_1c_2\(\underline m_1(z)-\underline m_2(z)\)\int\frac {t^2}{\(1+\underline m_1(z)t\)\(1+\underline m_2(z)t\)}dH(t)}{\left[c_1\int\frac t{1+\underline m_1(z)t}dH(t)-z\right]\left[c_1\int\frac t{1+\underline m_2(z)t}dH(t)-z\right]}.
\end{align*}
If $\underline m_1(z) \neq \underline m_2(z)$, we have
$$
\begin{aligned}
	 &\frac{c_1c_2\int\frac {t^2}{\(1+\underline m_1(z)t\)\(1+\underline m_2(z)t\)}dH(t)}{\left[c_1\int\frac t{1+\underline m_1(z)t}dH(t)-z\right]\left[c_1\int\frac t{1+\underline m_2(z)t}dH(t)-z\right]}=1.
\end{aligned}
$$
This, by the Cauchy inequality, yields
\begin{align}\label{uni2}
 \left(\frac{c_1c_2\int\frac {t^2}{\left|1+\underline m_1(z)t\right|^2}dH(t)}{\left|c_1\int\frac t{1+\underline m_1(z)t}dH(t)-z\right|^2}\cdot\frac{c_1c_2\int\frac {t^2}{\left|1+\underline m_2(z)t\right|^2}dH(t)}{\left|c_1\int\frac t{1+\underline m_2(z)t}dH(t)-z\right|^2} \right)^{\frac{1}{2}}\ge1.
\end{align}
On the other hand, we have from \eqref{uni1}
$$
\begin{aligned}
	\Im \underline m_j(z) & =\frac{\Im \underline m_j(z)c_1c_2\int\frac {t^2}{\left|1+\underline m_j(z)t\right|^2}dH(t)+v}{\left|c_1\int\frac t{1+\underline m_j(z)t}dH(t)-z\right|^2},
\end{aligned}
$$
which implies that for both $j=1$ and 2,
\begin{align}\label{uni3}
	\frac{c_1c_2\int\frac {t^2}{\left|1+\underline m_j(z)t\right|^2}dH(t)}{\left|c_1\int\frac t{1+\underline m_j(z)t}dH(t)-z\right|^2}<1.
\end{align}
The contradiction between \eqref{uni2} and \eqref{uni3} proves that $m_1(z) = m_2(z)$, and hence we are done.

\section{Proof of Theorems \ref{identitytestf} and Theorems \ref{identitytestl}}\label{cltappli}
We combine the proofs of Theorems \ref{identitytestf} and \ref{identitytestl} since both rely on applying the residue theorem and similar calculations.
Under the null hypothesis  $H_0$, it is known that
\begin{align*}
	z=-\frac{c_{n2}}{\underline m_n^0(z)}+\frac{c_{n1}}{1+\underline m_n^0(z)}.
\end{align*}
Let $\underline m_n^0(z)=-\frac1{1+\sqrt{\frac{c_{n1}}{c_{n2}}}r\xi}$, then we find 
\begin{align*}
	z=c_{n2}+\sqrt{c_{n1}c_{n2}}r\xi+\sqrt{c_{n1}c_{n2}}r^{-1}\xi^{-1}+c_{n1},\ {\rm and}\ dz=\sqrt{c_{n1}c_{n2}}\(r-r^{-1}\xi^{-2}\)d\xi,
\end{align*}
where $|\xi|=1$, and $r>1$ but close to $1$. When $\xi$ traverses the unit circle counterclockwise, $z$ will follow a counterclockwise path around a contour that encloses the support interval $\left[\left(\sqrt{c_{n2}}-\sqrt{c_{n1}}\right)^2,1\right]$, excluding the point $0$.

By performing calculations, we obtain the mean and covariance functions as follows
\begin{align*}
	\re X_f
	=&\lim_{r\downarrow1}\frac{1}{2\pi i\sqrt{c_{n1}c_{n2}}}\oint_{|\xi|=1}f(|\sqrt{c_{n2}}+\sqrt{c_{n1}}\xi|^2)
		\(\frac{\alpha_y }{\xi\(\xi^2-r^{-2}\)}+\frac{\Delta_y  }{\xi^3}\)d\xi\\
\triangleq&\frac{\alpha_y}{\sqrt{c_{n1}c_{n2}}}\mathcal E_1(f)+\frac{\Delta_y}{\sqrt{c_{n1}c_{n2}}}\mathcal E_2(f),
\end{align*}
and
\begin{align*}
	&{\rm Cov}\(X_{f_j},X_{f_k}\)\\
	=&-\frac{1+\alpha_y}{c_{n1}c_{n2}}\lim_{r\downarrow1}\frac{1}{4\pi^2}\oint_{|\xi_1|=1}\oint_{|\xi_2|=1}\frac{f_j\(|\sqrt{c_{n2}}+\sqrt{c_{n1}}\xi_1|^2\)f_k\(|\sqrt{c_{n2}}+\sqrt{c_{n1}}\xi_2|^2\)}{\(\xi_1-r\xi_2\)^2}d\xi_1d\xi_2\\
	&+\frac{\Delta_y}{c_{n1}c_{n2}}\oint_{|\xi_1|=1}\frac{f_j\(|\sqrt{c_{n2}}+\sqrt{c_{n1}}\xi_1|^2\)}{2\pi i\xi_1^2}d\xi_1\oint_{|\xi_2|=1}\frac{f_k\(|\sqrt{c_{n2}}+\sqrt{c_{n1}}\xi_2|^2\)}{2\pi i\xi_2^2}d\xi_2\\
	\triangleq&\frac{1+\alpha_y}{{c_{n1}c_{n2}}}\mathcal V_1(f_j,f_k)+\frac{\Delta_y}{{c_{n1}c_{n2}}}\mathcal V_2(f_j)\mathcal V_2(f_k).
\end{align*}
Let $f_1=x$, $f_2=x^2$, $f_3=\log(x).$ With these functions defined, we can now deduce
\begin{align*}
|\sqrt{c_{n2}}+\sqrt{c_{n1}}\xi|^2&={c_{n1}}+{c_{n2}}+\sqrt{c_{n1}c_{n2}}\(\xi+\xi^{-1}\),\\
|\sqrt{c_{n2}}+\sqrt{c_{n1}}\xi|^4&=\(c_{n1}+c_{n2}\)^2+2c_{n1}c_{n2}+2\sqrt{c_{n1}c_{n2}}\(c_{n1}+c_{n2}\)\(\xi+\xi^{-1}\)\\
&\qquad+c_{n1}c_{n2}\(\xi^2+\xi^{-2}\),\\
\log\(|\sqrt{c_{n2}}+\sqrt{c_{n1}}\xi|^2\)&=\log\(\sqrt{c_{n2}}+\sqrt{c_{n1}}\xi\)+\log\(\sqrt{c_{n2}}+\sqrt{c_{n1}}\xi^{-1}\).
\end{align*}
Through the application of the residue theorem, we derive that
\begin{align*}
	\mathcal E_1(f_1)=&\lim_{r\downarrow1}\frac{1}{2\pi i}\oint_{|\xi|=1} |\sqrt{c_{n2}}+\sqrt{c_{n1}}\xi|^2
	\(\frac{\xi}{\xi^2-r^{-2}}-\frac1{\xi}\)d\xi\\
=&\lim_{r\downarrow1}\frac{1}{2\pi i}\oint_{|\xi|=1}
\frac{\xi\({c_{n1}}+{c_{n2}}\)}{\xi^2-r^{-2}}d\xi+\lim_{r\downarrow1}\frac{1}{2\pi i}\oint_{|\xi|=1}
\frac{\xi^2\sqrt{c_{n1}c_{n2}}}{\xi^2-r^{-2}}d\xi\\	
&+\lim_{r\downarrow1}\frac{1}{2\pi i}\oint_{|\xi|=1}
\frac{\sqrt{c_{n1}c_{n2}}}{\xi^2-r^{-2}}d\xi-\frac{1}{2\pi i}\oint_{|\xi|=1}
\frac{{c_{n1}}+{c_{n2}}}{\xi}d\xi=0,\\
\mathcal E_1(f_2)=&\lim_{r\downarrow1}\frac{1}{2\pi i}\oint_{|\xi|=1}
\frac{\xi\(\(c_{n1}+c_{n2}\)^2+2c_{n1}c_{n2}\)}{\xi^2-r^{-2}}d\xi\\
&+\lim_{r\downarrow1}\frac{1}{2\pi i}\oint_{|\xi|=1} 
\frac{c_{n1}c_{n2}\xi^3}{\xi^2-r^{-2}}d\xi-\frac{1}{2\pi i}\oint_{|\xi|=1}\frac{\(c_{n1}+c_{n2}\)^2+2c_{n1}c_{n2}}{\xi}d\xi\\
=&c_{n1}c_{n2},\\
\mathcal E_1(f_3)=&\lim_{r\downarrow1}\frac{1}{2\pi i}\oint_{|\xi|=1} \log\(\sqrt{c_{n2}}+\sqrt{c_{n1}}\xi\)
\(\frac{\xi}{\xi^2-r^{-2}}-\frac1{\xi}\)d\xi\\
&+\lim_{r\downarrow1}\frac{1}{2\pi i}\oint_{|\xi|=1} \log\(\sqrt{c_{n2}}+\sqrt{c_{n1}}\xi^{-1}\)
\(\frac{\xi}{\xi^2-r^{-2}}-\frac1{\xi}\)d\xi\\
=&\lim_{r\downarrow1}\frac{1}{2\pi i}\oint_{|\xi|=1} \log\(\sqrt{c_{n2}}+\sqrt{c_{n1}}\xi\)
\(\frac{\xi}{\xi^2-r^{-2}}-\frac1{\xi}\)d\xi\\
=&\frac12\log\({c_{n2}}-c_{n1}\)-\frac12\log\({c_{n2}}\)=\frac12\log\(1-\frac{c_{n1}}{c_{n2}}\).
\end{align*}
Moreover, one finds
\begin{align*}
\mathcal E_2(f_1)=&\frac{1}{2\pi i}\oint_{|\xi|=1}|\sqrt{c_{n2}}+\sqrt{c_{n1}}\xi|^2\frac{1}{\xi^3}d\xi=0,\\
\mathcal E_2(f_2)=&\frac{1}{2\pi i}\oint_{|\xi|=1}\frac{c_{n1}c_{n2}}{\xi}d\xi=c_{n1}c_{n2},\\
\mathcal E_2(f_3)=&\frac{1}{2\pi i}\oint_{|\xi|=1}\log\(\sqrt{c_{n2}}+\sqrt{c_{n1}}\xi\)\frac{1}{\xi^3}d\xi=-\frac{c_{n1}}{2c_{n2}}.
\end{align*}
Hence, we see
\begin{align*}
	\re X_{f_1}
	=&0,\qquad
	\re X_{f_2}
	=\sqrt{c_{n1}c_{n2}}\(\alpha_y+\Delta_y\),\\
\re X_{f_3}
	=&\frac{\alpha_y\log\(1-\frac{c_{n1}}{c_{n2}}\)-\Delta_y\frac{c_{n1}}{c_{n2}}}{2\sqrt{c_{n1}c_{n2}}}.
\end{align*}

We are now in position to derive ${\rm Cov}\(X_{f_j},X_{f_k}\),j,k=1,2,3$. Let
\begin{align*}
	\mathcal V_3(f_j,\xi_2)=\lim_{r\downarrow1}\frac{1}{2\pi i}\oint_{|\xi_1|=1}\frac{f_j\(|\sqrt{c_{n2}}+\sqrt{c_{n1}}\xi_1|^2\)}{\(\xi_1-r\xi_2\)^2}d\xi_1,
\end{align*}
then
\begin{align*}
	\mathcal V_1(f_j,f_k)=\frac{1}{2\pi i}\oint_{|\xi_2|=1}f_k\(|\sqrt{c_{n2}}+\sqrt{c_{n1}}\xi_2|^2\)\mathcal V_3(f_j,\xi_2)d\xi_2.
\end{align*}
It is obvious that for fixed $|\xi_2|=1$, $|r\xi_2|=|r|>1$. So $r\xi_2$ is not a pole. Utilizing this fact, it follows that
\begin{align*}
	\mathcal V_3(f_1,\xi_2)=&\lim_{r\downarrow1}\frac{1}{2\pi i}\oint_{|\xi_1|=1}\frac{\sqrt{c_{n1}c_{n2}}}{\xi_1\(\xi_1-r\xi_2\)^2}d\xi_1=\frac{\sqrt{c_{n1}c_{n2}}}{\xi_2^2},\\
	\mathcal V_3(f_2,\xi_2)=&\lim_{r\downarrow1}\frac{1}{2\pi i}\oint_{|\xi_1|=1}\left[\frac{2\sqrt{c_{n1}c_{n2}}\(c_{n1}+c_{n2}\)}{\xi_1\(\xi_1-r\xi_2\)^2}+\frac{c_{n1}c_{n2}}{\xi_1^2\(\xi_1-r\xi_2\)^2}\right]d\xi_1\\
	=&\frac{2\sqrt{c_{n1}c_{n2}}\(c_{n1}+c_{n2}\)}{\xi_2^2}+\frac{2c_{n1}c_{n2}}{\xi_2^3},\\
	\mathcal V_3(f_3,\xi_2)=&\lim_{r\downarrow1}\frac{1}{2\pi i}\oint_{|\xi_1|=1}\frac{\log\(\sqrt{c_{n2}}+\sqrt{c_{n1}}\xi_1^{-1}\)}{\(\xi_1-r\xi_2\)^2}d\xi_1\\
	=&\lim_{r\downarrow1}\frac{1}{2\pi i\xi_2^2}\oint_{|\xi|=1}\frac{\log\(\sqrt{c_{n2}}+\sqrt{c_{n1}}\xi\)}{\(\xi-r^{-1}\xi_2^{-1}\)^2}d\xi=\frac{\sqrt{c_{n1}}}{\xi_2\(\sqrt{c_{n2}}\xi_2+\sqrt{c_{n1}}\)}.
\end{align*}
These imply 
\begin{align*}
	\mathcal V_1(f_1,f_1)=&\frac{1}{2\pi i}\oint_{|\xi_2|=1}\frac{{c_{n1}c_{n2}}}{\xi_2}d\xi_2={c_{n1}c_{n2}},\\
	\mathcal V_1(f_1,f_2)=&\frac{1}{2\pi i}\oint_{|\xi_2|=1}\frac{2{c_{n1}c_{n2}}\(c_{n1}+c_{n2}\)}{\xi_2}d\xi_2=2{c_{n1}c_{n2}}\(c_{n1}+c_{n2}\),\\
	\mathcal V_1(f_1,f_3)=&\frac{1}{2\pi i}\oint_{|\xi_2|=1}\log\(\sqrt{c_{n2}}+\sqrt{c_{n1}}\xi_2\)\frac{\sqrt{c_{n1}c_{n2}}}{\xi_2^2}d\xi_2=c_{n1},\\
	\mathcal V_1(f_2,f_2)=&\frac{1}{2\pi i}\oint_{|\xi_2|=1}\left[\frac{4{c_{n1}c_{n2}}\(c_{n1}+c_{n2}\)^2}{\xi_2}+\frac{2c_{n1}^2c_{n2}^2}{\xi_2}\right]d\xi_2\\
	=&{4{c_{n1}c_{n2}}\(c_{n1}+c_{n2}\)^2}+{2c_{n1}^2c_{n2}^2},\\
	\mathcal V_1(f_3,f_3)=&\frac{1}{2\pi i}\oint_{|\xi_2|=1}\frac{\sqrt{c_{n1}}\log\(\sqrt{c_{n2}}+\sqrt{c_{n1}}\xi_2\)}{\xi_2\(\sqrt{c_{n2}}\xi_2+\sqrt{c_{n1}}\)}d\xi_2
	=-\log\(1-c_{n1}/c_{n2}\).
\end{align*}
Furthermore, we deduce that
\begin{align*}
\mathcal V_2(f_1)=
	&\frac1{2\pi i}\oint_{|\xi|=1}\frac{\sqrt{c_{n1}c_{n2}}}{\xi}d\xi=\sqrt{c_{n1}c_{n2}},\\
\mathcal V_2(f_2)=
	&\frac1{2\pi i}\oint_{|\xi|=1}\frac{2\sqrt{c_{n1}c_{n2}}\(c_{n1}+c_{n2}\)}{\xi}d\xi=2\sqrt{c_{n1}c_{n2}}\(c_{n1}+c_{n2}\),\\
\mathcal V_2(f_3)=
	&\frac1{2\pi i}\oint_{|\xi|=1}\frac{\log\(\sqrt{c_{n2}}+\sqrt{c_{n1}}\xi\)}{\xi^2}d\xi=\frac{\sqrt{c_{n1}}}{\sqrt{c_{n2}}}.
\end{align*}
Therefore, we get
\begin{align*}
{\rm Cov}\(X_{f_1},X_{f_1}\)
=&1+\alpha_y+\Delta_y,\\
{\rm Cov}\(X_{f_1},X_{f_2}\)
=&2\(c_{n1}+c_{n2}\)\(1+\alpha_y+\Delta_y\),\\
{\rm Cov}\(X_{f_1},X_{f_3}\)
=&\frac{1+\alpha_y}{c_{n2}}+\frac{\Delta_y}{c_{n2}},\\
{\rm Cov}\(X_{f_2},X_{f_2}\)
=&4\(c_{n1}+c_{n2}\)^2\(1+\alpha_y+\Delta_y\)+{2c_{n1}c_{n2}}\(1+\alpha_y\),\\
{\rm Cov}\(X_{f_3},X_{f_3}\)
=&-\frac{1+\alpha_y}{{c_{n1}c_{n2}}}\log\(1-c_{n1}/c_{n2}\)+\frac{\Delta_y}{c_{n2}^2}.
\end{align*}

Denote 
\begin{align*}\boldsymbol\zeta_{n}&=\begin{pmatrix}
	\zeta_{1n}\\
	\zeta_{2n}\\
	\zeta_{3n}
\end{pmatrix}=\begin{pmatrix}
	\frac{\nu_{np}}{\sqrt{pn}}\int xdF^{\bS_n}(x)\\
	\frac{\nu_{np}}{\sqrt{pn}}\int x^2dF^{\bS_n}(x)\\
	\frac{\nu_{np}}{\sqrt{pn}}\int \log(x)dF^{\bS_n}(x)
\end{pmatrix},
\boldsymbol{\theta}=\begin{pmatrix}
	\theta_{1}\\
	\theta_{2}\\
	\theta_3
\end{pmatrix}=\begin{pmatrix}
\frac{\nu_{np}}{\sqrt{pn}}\int xdF^{c_{n1},c_{n2},H_n}(x)\\
\frac{\nu_{np}}{\sqrt{pn}}\int x^2dF^{c_{n1},c_{n2},H_n}(x)\\
\frac{\nu_{np}}{\sqrt{pn}}\int \log(x)dF^{c_{n1},c_{n2},H_n}(x)
\end{pmatrix}
.\end{align*}
Subsequently, we proceed to compute
\begin{align*}
	\theta_1=&\frac{\nu_{np}}{\sqrt{pn}}\frac1{2\pi i}\oint\frac{1}{1+\underline m_n^0(z)}dz=\sqrt{\frac{c_{n2}}{c_{n1}}}\frac1{2\pi i}\oint_{|\xi|=1}\frac{1}{\xi}d\xi=\sqrt{\frac{c_{n2}}{c_{n1}}},\\
	\theta_2=&\frac{\nu_{np}}{\sqrt{pn}}\frac1{2\pi i}\oint\frac{z}{1+\underline m_n^0(z)}dz\\
		=&\frac1{2\pi i}\oint_{|\xi|=1}\frac{\(\xi^{2}-1\)|\sqrt{c_{n2}}+\sqrt{c_{n1}}\xi|^2}{\xi^2}d\xi\\
		&+\sqrt{\frac{c_{n2}}{c_{n1}}}\frac1{2\pi i}\oint_{|\xi|=1}\frac{\(\xi^{2}-1\)|\sqrt{c_{n2}}+\sqrt{c_{n1}}\xi|^2}{\xi^3}d\xi\\
		=&\sqrt{\frac{c_{n2}}{c_{n1}}}\frac1{2\pi i}\oint_{|\xi|=1}\frac{{c_{n2}}+{c_{n1}}}{\xi}d\xi=\sqrt{\frac{c_{n2}}{c_{n1}}}\(c_{n1}+c_{n2}\),
\end{align*}
and
\begin{align*}
	\theta_3=&\frac{\nu_{np}}{\sqrt{pn}}\frac1{2\pi i}\oint\frac{\log(z)}{z\(1+\underline m_n^0(z)\)}dz\\
	=&\sqrt{\frac{c_{n2}}{c_{n1}}}\frac1{2\pi i}\oint_{|\xi|=1}\frac{\(\xi^{2}-1\)\log\(|\sqrt{c_{n2}}+\sqrt{c_{n1}}\xi|^2\)}{\(\sqrt{c_{n2}}+\sqrt{c_{n1}}\xi\)\(\sqrt{c_{n2}}\xi+\sqrt{c_{n1}}\)\xi^2}d\xi\\
	&+\frac1{2\pi i}\oint_{|\xi|=1}\frac{\(\xi^{2}-1\)\log\(|\sqrt{c_{n2}}+\sqrt{c_{n1}}\xi|^2\)}{\(\sqrt{c_{n2}}+\sqrt{c_{n1}}\xi\)\(\sqrt{c_{n2}}\xi+\sqrt{c_{n1}}\)\xi}d\xi\\
	=&\sqrt{\frac{c_{n2}}{c_{n1}}}\frac1{\pi i}\oint_{|\xi|=1}\frac{\log\(\sqrt{c_{n2}}+\sqrt{c_{n1}}\xi\)}{\(\sqrt{c_{n2}}+\sqrt{c_{n1}}\xi\)\(\sqrt{c_{n2}}\xi+\sqrt{c_{n1}}\)}d\xi\\
	&-\sqrt{\frac{c_{n2}}{c_{n1}}}\frac1{2\pi i}\oint_{|\xi|=1}\frac{\log\(\sqrt{c_{n2}}+\sqrt{c_{n1}}\xi\)}{\xi^2\(\sqrt{c_{n2}}+\sqrt{c_{n1}}\xi\)\(\sqrt{c_{n2}}\xi+\sqrt{c_{n1}}\)}d\xi\\
	&-\sqrt{\frac{c_{n2}}{c_{n1}}}\frac1{2\pi i}\oint_{|z|=1}\frac{z^2\log\(\sqrt{c_{n2}}+\sqrt{c_{n1}}z\)}{\(\sqrt{c_{n2}}z+\sqrt{c_{n1}}\)\(\sqrt{c_{n2}}+\sqrt{c_{n1}}z\)}dz\\
	=&\sqrt{\frac{c_{n2}}{c_{n1}}}\frac{2\log\({c_{n2}}-c_{n1}\)-\log\({c_{n2}}\)}{\({c_{n2}}-{c_{n1}}\)}\\
	&-\({\frac{c_{n2}}{c_{n1}}}\)^{3/2}\frac{\log\({c_{n2}}-c_{n1}\)-\frac12\log\({c_{n2}}\)}{{c_{n2}}-{c_{n1}}}\\
	&-\sqrt{\frac{c_{n2}}{c_{n1}}}\frac{c_{n1}-\frac12\({c_{n1}}+{c_{n2}}\)\log\({c_{n2}}\)}{{c_{n1}}{c_{n2}}}\\
	&-\sqrt{\frac{c_{n1}}{c_{n2}}}\frac{\log\({c_{n2}}-c_{n1}\)-\frac12\log\({c_{n2}}\)}{\({c_{n2}}-{c_{n1}}\)}\\
	=
	&\frac{1}{\sqrt{c_{n1}c_{n2}}}\(\log\({c_{n2}}-c_{n1}\)-1\)+\sqrt{\frac{c_{n2}}{c_{n1}}}\frac{\log\({c_{n2}}\)-\log\({c_{n2}}-c_{n1}\)}{c_{n1}}.
\end{align*}
From Theorem \ref{cltthm}, we get the following results
\begin{align*}
	&p\begin{pmatrix}
		\zeta_{1n}-\theta_1\\
		\zeta_{2n}-\theta_2
	\end{pmatrix}\xrightarrow{d}N\Bigg[\begin{pmatrix}
		0\\
		\sqrt{c_{n1}c_{n2}}\(\alpha_y+\Delta_y\)
	\end{pmatrix},\begin{pmatrix}
		{\rm Cov}\(X_{f_1},X_{f_1}\)&{\rm Cov}\(X_{f_1},X_{f_2}\)\\
		{\rm Cov}\(X_{f_1},X_{f_2}\)& {\rm Cov}\(X_{f_2},X_{f_2}\)
	\end{pmatrix}\Bigg]
\end{align*}
and
\begin{align*}
	&p\begin{pmatrix}
		\zeta_{1n}-\theta_1\\
		\zeta_{3n}-\theta_3
	\end{pmatrix}\xrightarrow{d}N\Bigg[\begin{pmatrix}
		0\\
		\frac{\alpha_y\log\(1-\frac{c_{n1}}{c_{n2}}\)-\Delta_y\frac{c_{n1}}{c_{n2}}}{2\sqrt{c_{n1}c_{n2}}}
	\end{pmatrix},\begin{pmatrix}
		{\rm Cov}\(X_{f_1},X_{f_1}\)&{\rm Cov}\(X_{f_1},X_{f3}\)\\
		{\rm Cov}\(X_{f_1},X_{f_3}\)& {\rm Cov}\(X_{f_3},X_{f_3}\)
	\end{pmatrix}\Bigg].
\end{align*}
Notice that the statistics $nW$ and $\mathcal L$ can be represented as
\begin{align*}
	nW=p\left[\frac1{\sqrt{c_{n1}c_{n2}}}\zeta_{2n}-2\sqrt{\frac{c_{n2}}{c_{n1}}}\zeta_{1n}-\frac{c_{n1}}{c_{n2}}\zeta_{1n}^2+1+\frac{c_{n2}}{c_{n1}}\right]
\end{align*}
and
\begin{align*}
	\mathcal L=p\left[\sqrt{\frac{c_{n1}}{c_{n2}}}\zeta_{1n}-\sqrt{c_{n1}c_{n2}}\zeta_3+\log\(c_{n2}\)-1\right].
\end{align*}
Let
\begin{align*}
	g_1\(\boldsymbol{\theta}\)=\frac1{\sqrt{c_{n1}c_{n2}}}\theta_{2}-2\sqrt{\frac{c_{n2}}{c_{n1}}}\theta_{1}-\frac{c_{n1}}{c_{n2}}\theta_{1}^2+1+\frac{c_{n2}}{c_{n1}},
\end{align*}
and
\begin{align*}
	g_2\(\boldsymbol{\theta}\)=\sqrt{\frac{c_{n1}}{c_{n2}}}\theta_1-\sqrt{c_{n1}c_{n2}}\theta_3+\log\(c_{n2}\)-1.
\end{align*}
Applying the Delta method, we deduce
\begin{align*}
	nW-p=p\(g_1\(\boldsymbol{\zeta}\)-g_1\(\boldsymbol{\theta}\)\)\xrightarrow{d}N\(\(\alpha_y+\Delta_y\),2\(1+\alpha_y\)\),
\end{align*}
and
\begin{align*}
p\(g_2\(\boldsymbol{\zeta}\)-g_2\(\boldsymbol{\theta}\)\)\xrightarrow{d}N\(-\frac{\alpha_y}{2}\log\(1-\frac{c_{n1}}{c_{n2}}\)+\frac{\Delta_yc_{n1}}{2c_{n2}},-\(1+\alpha_y\)\(\frac{c_{n1}}{c_{n2}}+\log\(1-\frac{c_{n1}}{c_{n2}}\)\)\).
\end{align*}
Here
\begin{align*}
	g_2(\boldsymbol\theta)=1-\(1-\frac{c_{n2}}{c_{n1}}\)\log\(1-\frac{c_{n1}}{c_{n2}}\).
\end{align*}
Thus, we conclude the proof of these two theorems.

\section{Verification of Remark \ref{rem1}}\label{verirem2}
To verify Remark \ref{rem1}, we must truncate $y_{jk}$ at a positive constant $C$ and subsequently renormalize it. This procedure will lead to the following result
\begin{align*}
\re\left|\bx_1^*\bB\bx_1-\rtr\(\bB\bSig_n\)\right|^q
	\le& C_q\left[\(\rtr\(\bB\bSig_n\bB^*\bSig_n\)\)^{q/2}+\rtr\(\bB\bSig_n\bB^*\bSig_n\)^{q/2}\right]
	=O\(p^{q/2}\).
\end{align*}

\subsection{Truncation}
Let $\hat {y}_{jk}=y_{jk}I\left(\left|y_{jk}\right|\leq C\right)$, $\widehat\bY_n=\left(\hat {y}_{jk}\right)$ and
$$\widehat{\mathbf S}_n=\frac{1}{\nu_{np}}\bGa_n{\widehat {\mathbf Y}_n}{\widehat {\mathbf Y}_n^*}\bGa_n^*.$$
By employing Lemma \ref{lemma:4}, one can derive
	\begin{align*}
	{L^4}\left({F^{\mathbf S_n}},{F^{\widehat{\mathbf S}_n}}\right) \le& \frac{2\|\bSig_n\|^2}{{{p^2\nu_{np}^2}}}{\rm tr}\left(\bY_n\bY_n^*+\widehat\bY_n\widehat\bY_n^*\right){\rm tr}\left[\left(\bY_n-\widehat\bY_n\right){\left(\bY_n-\widehat\bY_n\right)^ * }\right]\\
\le&\frac{4K^2}{{{p^2n^2}}}\left(\sum_{jk}|y_{jk}|^2\right)\left(\sum_{jk}|y_{jk}|^2I\(|y_{jk}|>C\)\right)\\
\xrightarrow{a.s.}&{4K^2}\re|y_{11}|^2I\(|y_{11}|>C\)	
	\end{align*}
where $\|\bSig_n\|\le K$. By choosing $C$ large enough, $\re|y_{11}|^2I\(|y_{11}|>C\)$ can be made arbitrarily small.
\subsection{Centralization}
	Denote $$\widetilde {y}_{jk}=\widehat {y}_{jk}-{\rm E}\(\widehat {y}_{jk}\), \ \widetilde\bY_n=\left(\widetilde {y}_{jk}\right)~~\mbox{and}~~ \ \widetilde{\mathbf S}_n=\frac{1}{\nu_{np}}\bGa_n\widetilde\bY_n\widetilde\bY_n^*\bGa_n^*.$$
Using Lemma \ref{lemma:3}, we have
\begin{align*}
		\left\| F^{\widehat\bS_n} -F^{\widetilde\bS_n}\right\|_{KS} \le& \frac{1}{p}{\rm rank}\left(\frac{1}{\sqrt{\nu_{np}}}\bGa_n\widehat\bY_n-\frac{1}{\sqrt{\nu_{np}}}\bGa_n\widetilde\bY_n\right)\\
		\le& \frac{1}{p}{\rm rank}\left(\widehat\bY_n-\widetilde\bY_n\right)=\frac{1}{p}{\rm rank}\left(\re\widehat\bY_n\right)\to0.
\end{align*}
\subsection{Rescaling}
Define $\sigma^2={\rm E}\left|\tilde y_{11}\right|^2.$
Owing to Lemma \ref{lemma:4}, we get
	\begin{align*}
		{L^4}\left(F^{\widetilde {\mathbf S}_n} , F^{\sigma^{-2}\widetilde {\mathbf S}_n}\right)
		\le &\frac{2K^2\(1+\sigma^{-2}\)\(\sigma^{-2}-1\)}{p^2\nu_{np}^2}{\rm tr}^2\left(\widetilde{\bY}_n\widetilde{\bY}_n^*\right)\notag\\
		\le &\frac{2K^2\(1+\sigma^{-2}\)\(\sigma^{-2}-1\)}{p^2n^2}\left(\sum_{jk}|\tilde y_{jk}|^2\right)^2\notag\\
		\xrightarrow{a.s.}&2K^2\(1-\sigma^{4}\).
	\end{align*}
Because $\sigma^2\to1$ as $C\to\infty$, we know that $1-\sigma^{4}$ can be made arbitrarily small by choosing $C$ large enough.


\end{document}